\documentclass[11pt]{amsart}
\usepackage{amsxtra}

\theoremstyle{plain}

\def\CC {{\mathbb C}}
\def\RR {{\mathbb R}}
\def\NN {{\mathbb N}}
\def\ZZ {{\mathbb Z}}
\def\PP {{\mathbb P}}

\def\be {\begin{eqnarray}}
\def\ben {\begin{eqnarray*}}
\def\ee {\end{eqnarray}}
\def\een {\end{eqnarray*}}

\def\AAA{\kern-0.3em}
\def\AA{\kern-0.18em}
\def\AC{\kern-0.14em}
\def\AB{\kern-0.22em}

\newcommand \nc {\newcommand}

\newtheorem{theorem}{Theorem}[section]
\newtheorem{lemma}[theorem]{Lemma}
\newtheorem{proposition}[theorem]{Proposition}
\newtheorem{corollary}[theorem]{Corollary}
\newtheorem{definition}[theorem]{Definition}
\newtheorem{example}[theorem]{Example}
\newtheorem{remark}[theorem]{Remark}
\newtheorem{conjecture}[theorem]{Conjecture}
\nc \bth[1] { \begin{theorem}\label{t#1} } \nc \ble[1] {
\begin{lemma}\label{l#1} } \nc \bpr[1] {
\begin{proposition}\label{p#1} } \nc \bco[1] {
\begin{corollary}\label{c#1} } \nc \bde[1] {
\begin{definition}\label{d#1}\rm } \nc \bex[1] {
\begin{example}\label{e#1}\rm } \nc \bre[1] {
\begin{remark}\label{r#1}\rm } \nc \bcon[1] { 
\begin{conjecture}\label{co#1}\rm} \nc \bqu[1]  {
\medskip\noindent{\it{Question #1}} }
\nc {\ethe} { \end{theorem} }
 \nc {\ele} { \end{lemma} } \nc {\epr}
{ \end{proposition} } \nc {\eco} { \end{corollary} } \nc {\ede} {
\end{definition} } \nc {\eex} { \end{example} } \nc {\ere} {
\end{remark} } \nc {\econ} { \end{conjecture} } \nc {\equ} {\smallskip}
 \nc \thref[1]{Theorem \ref{t#1}}
\nc \leref[1]{Lemma \ref{l#1}} \nc \prref[1]{Proposition
\ref{p#1}} \nc \coref[1]{Corollary \ref{c#1}} \nc
\deref[1]{Definition \ref{d#1}} \nc \exref[1]{Example \ref{e#1}}
\nc \reref[1]{Remark \ref{r#1}}
\nc \conref[1]{Conjecture \ref{co#1}}

\def \d {{\mathrm d}}

\def \B {{\mathcal B}}

\def \L {{\mathcal L}}
\def \R {{\mathcal Re}}

\def \diag { {\mathrm{diag}} }

 \def\AA  {\kern-0.1em}
 \def\BB  {\kern+0.1em}
 \def\BBB {\kern+0.15em}
 \def\K   {\kern+0.05em}
 \def\MK  {\kern-0.07em}
 \def\MKK {\kern-0.04em}

\begin{document}

\vspace{0.5cm}

\title[Stokes matrices]
{ Stokes matrices for a class of reducible equations }

\author[Tsvetana  Stoyanova]{Tsvetana  Stoyanova}

\date{24.06.2019}

 \maketitle

\begin{center}
{Department of Mathematics and Informatics,
Sofia University,\\ 5 J. Bourchier Blvd., Sofia 1164, Bulgaria, 
cveti@fmi.uni-sofia.bg}
\end{center}

\vspace{0.5cm}

{\bf Abstract.} This paper is a continuation of our previous work \cite{St}
where we have studied the Stokes phenomenon for a particular family of 
equation \eqref{initial} with  \eqref{form-0}-\eqref{npe} from a 
perturbative point of view. Here we focis on the explicit computation
of the Stokes matrices at the non-resonant irregular singularity
for a more general situation. In particular, utilizing Borel-Laplace summation method,
the iterated integrals approach and some
properties of the hypergeometric series  we compute by hand the Stokes matrices
of three  families of equation \eqref{initial}-\eqref{form-0}-\eqref{npe} 
under assumptions that $\beta_j$'s are assumed distinct
and $|\beta_3-\beta_1| < |\beta_3-\beta_2|$. Moreover these results remain valid
for these distinct $\beta_j$'s for which $|\beta_3-\beta_1|=|\beta_3-\beta_2|$
but $\beta_3-\beta_1 \neq \pm (\beta_3-\beta_2)$ on condition that
$\R  (\alpha_2-\alpha_1) > -1$. 
In addition, iterated integrals approach allows us to give, under some restrictions, an explicit representation
of the 1-sum of the product of two certain divergent 1-summable power series, that have
different singular directions.

{\bf Key words: Stokes matrices, Borel-Laplace summation, Iterated integrals, Hypergeometric series, 
	Third-order reducible  linear ordinary differential equation}

{\bf 2010 Mathematics Subject Classification:  34M40, 34M30, 34M03, 34A25}

\headsep 10mm \oddsidemargin 0in \evensidemargin 0in

\section{Introduction}

    We consider a linear reducible ordinary differential equation  
	  \be\label{initial}
		   L\,y=0,
			\ee
			in a complex variable $x$,
			where $L$ is a third-order linear differential operator of the form
			\be\label{form-0}
			  L=L_3 \circ L_2 \circ L_1, \quad
			  L_i \circ L_j \neq L_j \circ L_i.
			\ee
		The first order differential operators $L_j$ have the form	
			\be\label{npe}
		  L_j =\partial-\left(\frac{\alpha_j}{x} + \frac{\beta_j}{x^2}\right),\quad
		   \partial=\frac{d}{d x},
			\ee
			where $\alpha_j,\,\beta_j\in\CC$ such that
			$\beta_j$'s are assumed distinct. Throughout this paper we also assume
			that $\arg (\beta_1-\beta_2) \neq \arg (\beta_1-\beta_3) \neq \arg (\beta_2-\beta_3)
			\neq \arg (\beta_1-\beta_2)$.

			The equation \eqref{initial} with \eqref{form-0}-\eqref{npe} (in short, the equation \eqref{initial}) 
			is a third-order solvable differential equation, in the sense that its differential Galois
		group is a solvable linear algebraic group \cite{K}. This  equation, in general, has over
		 $\CC\PP^1$  two singular points - an irregular singular point at the origin
	  of Poincar\'e rank 1 and a regular singular point at $x=\infty$. 
      In \cite{St} we considered a particular family of equation \eqref{initial}
      for $\alpha_1=0, \alpha_2=\nu-2, \alpha_3=\nu-4, \beta_1=1, \beta_2=2, \beta_3=0$,
      where $\nu\in\CC$ is an arbitrary parameter. There we investigated  the Stokes
      phenomenon  applying a perturbative approach. With the present paper we continue our
      study on the Stokes phenomenon. The main goal of this paper is the explicit computation of the Stokes matrices
      of equation \eqref{initial} when the origin is a non-resonant irregular singularity.

     The Stokes matrices are analytic invariants of equation \eqref{initial}.
     They measure the difference between two different 1-sums of the same formal fundamental
     matrix $\hat{\Phi}(x)$ at the origin on each side of the so called singular directions.		
     There are many ways for the computation of the Stokes matrices - analytic, algebraic,
     topological - depending mainly on the kind of the considered equations. In this paper
     we utilize Borel-Laplace summation method. Generally speaking, starting from a formal matrix solution
     $\hat{\Phi}(x)$ at the origin we associate with it 
     an actual matrix solution $\Phi(x)$ on a sector at $x=0$ applying Borel-Laplace 
     summation approach  (see next section for details). 
     On this sector the actual matrix solution  $\Phi(x)$ represent an unique actual function,
     asymptotic at the origin to the formal matrix $\hat{\Phi}(x)$ in Gevrey order 1 sense. 
     Let $\theta$ be a singular direction of
     equation \eqref{initial}. Let $\Phi^{-}_{\theta}(x)=\Phi_{\theta-\epsilon}(x)$
     and $\Phi^{+}_{\theta}=\Phi_{\theta+\epsilon}(x)$ where $\epsilon$ is a small positive
     number, be two actual matrix solutions at the origin of equation \eqref{initial}
     obtained from $\hat{\Phi}(x)$ by the Borel-Laplace summation method. 
     The solution $\Phi^{-}_{\theta}(x)$ can not be extended analitycally to the solution
     $\Phi^{+}_{\theta}(x)$. It must jump to the solution $\Phi^{+}_{\theta}(x)$ and the
     Stokes matrix $St_{\theta}\in GL_3(\CC)$ measures this jump by
      $$\,
       \Phi^{-}_{\theta}(x)=\Phi^{+}_{\theta}(x)\,St_{\theta}.
      \,$$		
      Moreover, the matrix solution $\Phi^{+}_{\theta}(x)\,St_{\theta}$ is no longer
      asymptotic at the origin to the formal solution $\hat{\Phi}(x)$. This phenomenon
      is called Stokes phenomenon named for the first mathematician observed it.

    The problem of explicit computation of the Stokes matrices for a higer order scalar ordinary differential
    equation has not been attacked in so many works. In \cite{M-dP, dP}
    van der Put and Cruz Morales compute  the Stokes matrices for quantum differential equations using
    multisummation and monodromy identity. In \cite{DM, M} Duval and Mitschi deal with the Stokes matrices of confluent generalized
	 hypergeometric equations using analytic continuation of $G$-function of Meijer.
	In \cite{H} Hein computes the Stokes matrices of unramified confluent hypergeometric equations
	using a toplogical approach.
	In this paper we employ the Borel-Laplace summation approach in combination with
	some properties of the hypergeometric series for the computation of the Stokes
	matrices of equation \eqref{initial}. For a general irregular singularity of Poincar\'{e} rank 1
	the explicit computation of the Stokes matrices is a difficult task. The Borel-Laplace
	summation method is an useful tool to solve this problem in the case when the Borel
	transform of the formal solution $\hat{\Phi}(x)$ (in fact the Borel transform 
	of the folmal matrix $\hat{H}(x)$ such that $\hat{\Phi}(x)=\hat{H}(x)\,x^{\Lambda}\,\exp(-Q/x)$)
	 can be expressed in terms of the 
	special functions. But such differential equations for which it is possible are very special.
	Fortunately equation \eqref{initial} is rather special. The point is that such a reducible equation
	admits a global fundamental matrix whose off-diagonal elements are expressed in terms of iterated integrals.
	In this article we fully exploit such a basis in order to compute the Stokes matrices.
	As a result we provide explicit formulas for the Stokes multipliers for all values of
	the parameters $\alpha_j$ when $|\beta_3-\beta_1| < |\beta_3-\beta_2|$ and distinct $\beta_j$ except 
	the case when $\alpha_3-\alpha_1\notin\ZZ_{\leq -4}$ but $\alpha_3-\alpha_2\in\ZZ_{\geq -1}$
	such that $\alpha_2-\alpha_1\notin\ZZ_{\leq -2}$. 
	  This case  is out of scope of this work by techical reasons
	  (see \reref{transf}) and it is postponed for forther researches. In fact the exclusion of
	  this case of our work emphasizes that the explicit expression of the Borel transform
	  of $\hat{H}(x)$ in terms of classical functions is out of reach. 
	  Moreover, it turns out that the obtained formulas remain valid for all values of
	  the distinct $\beta$'s for which $|\beta_3-\beta_1|=|\beta_3-\beta_2|$
	  provided that $\beta_3-\beta_1 \neq \pm (\beta_3-\beta_2)$ and $\R (\alpha_2-\alpha_1) > -1$. 
	  Aside of these results
	  the iterated integrals approach also provides an useful tool for the computation by hand
	  of the 1-sum of the product of two divergent power series that have distinct singular
	  directions. 
	
     Summarily,	the basic results of this paper are as follows:\\
      In \thref{Stokes-2} we compute explicitly, by hand, the Stokes matrices at the
      non-resonant irregular singularity  $x=0$ when $\beta_j$'s are distinct
       such that $|\beta_3-\beta_1| < |\beta_3-\beta_2|$ and when
       \begin{itemize}
       	\item\,
       	$\alpha_3-\alpha_1\notin\ZZ_{\leq -4}$ but $\alpha_3-\alpha_2\notin\ZZ_{\geq -1}$.
       	
       	\item\,
       	$\alpha_3-\alpha_1\notin\ZZ_{\leq -4}$ and $\alpha_3-\alpha_2\in\ZZ_{\geq -1}$ but
       	$\alpha_2-\alpha_1\in\ZZ_{\leq -2}$.
       	
       	\item\,
       	$\alpha_3-\alpha_1\in\ZZ_{\leq -4}$.
       \end{itemize}
 In \thref{Stokes-3} we show that the results of \thref{Stokes-2} remain valid for all distinct $\beta_j$'s
 such that $|\beta_3-\beta_1|=|\beta_3-\beta_2|$ provided that $\beta_3-\beta_1 \neq \pm (\beta_3-\beta_2)$
 and $\R (\alpha_2-\alpha_1) > -1$ and when
     \begin{itemize}
     	\item\,
     	  	$\alpha_3-\alpha_1\notin\ZZ_{\leq -4}$ but $\alpha_3-\alpha_2\notin\ZZ_{\geq -1}$.
     	  	
     	\item\,
     	  	     	$\alpha_3-\alpha_1\in\ZZ_{\leq -4}$.
     \end{itemize}
 
     As a by-product we also compute explicitly, by hand, 
     the 1-sum of the product
      $$\,
      \hat{\omega}(x)=
      \left(\sum_{n=0}^{\infty} (-1)^n
      \frac{(2+\alpha_2-\alpha_1)^{(n)}}{(\beta_2-\beta_1)^n}\,x^n\right)
         \left(\sum_{n=0}^{\infty} (-1)^n
         \frac{(2+\alpha_3-\alpha_2)^{(n)}}{(\beta_3-\beta_2)^n}\,x^n\right)
      \,$$
      when $|\beta_3-\beta_1| < |\beta_3-\beta_2|,\,|\beta_3-\beta_1| < |\beta_2-\beta_1|$
      such that $\beta_j$'s are distinct and when
      \begin{itemize}
      	\item\,
      	$\alpha_3-\alpha_1\notin\ZZ_{\leq -4}$ such that $\alpha_2-\alpha_1,\,\alpha_3-\alpha_2\notin\ZZ$
      	(see \thref{sum}).
      	
      	\item\,
      	$\alpha_3-\alpha_1\in\ZZ_{\leq -4}$ but $\alpha_2-\alpha_1,\,\alpha_3-\alpha_2\notin\ZZ_{\leq -2}$
      	(see \thref{sum-1}).
      \end{itemize}
    \thref{sum} shows in a direct way that the 1-sum of the product $\hat{\omega}(x)$ has three different
    singular directions $\theta_1=\arg (\beta_1-\beta_2),\,\theta_2=\arg (\beta_2-\beta_3)$
    and $\theta_3=\arg (\beta_1-\beta_3)$ provided that $\theta_1 \neq \theta_2 \neq \theta_3 \neq \theta_1$. 
    Observe that $\theta_1$ and $\theta_2$ are the
    singular directions of the initial power series. The direction $\theta_3$ is an additional singular
    direction. This result exemplifies the statements of Loday-Richaud and Remy \cite{LR-R},
    and Sauzin \cite{Sa-1} that the set $\widehat{\mathcal{R}es}_{\Omega}$ of resurgent functions
    with singular support $\Omega$ is not stable under convolution.
    Moreover, since the iterated integrals naturally govern such products of power series,
    we observe the same phenomenon in the construction of the 1-sum $\Phi(x)$. More precisely,
    when  	$\alpha_3-\alpha_1\notin\ZZ_{\leq -4}$ but $\alpha_3-\alpha_2\notin\ZZ_{\geq -1}$
    and $\theta_2 \neq \arg (\beta_1-\beta_3)$ the element $\Phi_{13}(x)$ of $\Phi(x)$ 
    has the form (see \prref{as-2})
     $$\,
            \frac{x^{\alpha_3+4}e^{-\frac{\beta_3}{x}}}{(\beta_3-\beta_2)(\beta_2-\beta_1)}
            	\left[
            F(\alpha; \beta; 0)\,
            \int_0^{+\infty e^{i \theta_2}}
            \left(1+\frac{\nu}{\beta_3-\beta_1}\right)^{\alpha_1-\alpha_3-4}
            e^{-\frac{\nu}{x}} d \left(\frac{\nu}{x}\right)
            - \psi_{\theta}(x)\right],       
     \,$$
  where
  $$\,
    	F(\alpha; \beta; 0)=\sum_{p=0}^{\infty}
    	\frac{(2+\alpha_3-\alpha_2)^{(p)}}{(4+\alpha_3-\alpha_1)^{(p)}}
    	\left(\frac{\beta_3-\beta_1}{\beta_3-\beta_2}\right)^p    
  \,$$
 with $(a)^{(n)}=a(a+1) \dots (a+n-1),\,(a)^{(0)}=1$. The function $\psi_{\theta}(x)$ is given by
    	\ben
    	& &
    	\psi_{\theta}(x) =
    	- \frac{\beta_3-\beta_1}{\beta_3-\beta_2}
    	\frac{\Gamma(4+\alpha_3-\alpha_1) \Gamma(\alpha_2-\alpha_3-1)}
    	{\Gamma(2+\alpha_2-\alpha_1)}
    	\left(\frac{\beta_1-\beta_2}{\beta_3-\beta_2}\right)^{\alpha_2-\alpha_3-3}\\[0.4ex]
    	&\times&
    	\left(\frac{\beta_1-\beta_3}{\beta_1-\beta_2}\right)^{\alpha_1-\alpha_3-4}
    	\int_0^{+\infty e^{i \theta_2}}
    	\left(1+\frac{\xi}{\beta_3-\beta_1}\right)^{\alpha_1-\alpha_3-4}
    	e^{-\frac{\xi}{x}} d\left(\frac{\xi}{x}\right)\\[0.85ex]  
    	&+&
    	\frac{\beta_3-\beta_1}{\beta_1-\beta_2}
    	\sum_{s=0}^{\infty}
    	\frac{(2+\alpha_2-\alpha_1)^{(s)}}{(\alpha_2-\alpha_3-1)^{(s)}}
    	\left(\frac{\beta_3-\beta_2}{\beta_1-\beta_2}\right)^s
    	\int_0^{+\infty\,e^{i \theta_3}}
    	\left(1+\frac{\xi}{\beta_3-\beta_2}\right)^{\alpha+s}\,
    	e^{-\frac{\xi}{x}}\,d \left(\frac{\xi}{x}\right),
    	\een
    	where $\alpha=\alpha_2-\alpha_3-2$, $\Gamma(z)$ is Euler's Gamma function
    	and $\theta_3 \neq \arg (\beta_2-\beta_3)$.
  The first term of $\Phi_{13}(x)$ and its singular direction $\theta=\arg (\beta_1-\beta_3)$
  are predetermined from the equation. But we are rather surprise to see  the same singular direction in
  the function $\psi_{\theta}(x)$. Its existence comes from the relation between $\psi_{\theta}(x)$
  and the 1-sum of the  product $\hat{\omega}{(x)}$. In particular, the function $\psi_{\theta}(x)$ is the 1-sum of 
  the series
   $$\,
    \hat{\psi}(x)=\sum_{n=0}^{\infty} (-1)^n
    \frac{(2+\alpha_3-\alpha_2)^{(n)}}{(\beta_3-\beta_2)^n}
    \left(\sum_{p=0}^{\infty}
     \frac{(2+n+\alpha_3-\alpha_2)^{(p)}}{(4+n+\alpha_3-\alpha_1)^{(p)}}
     \left(\frac{\beta_3-\beta_1}{\beta_3-\beta_2}\right)^p\right)\,x^n
   \,$$
  when $|\beta_3-\beta_1| < |\beta_3-\beta_2|,\,\alpha_3-\alpha_1\notin\ZZ_{\leq -4}$
  but $\alpha_3-\alpha_2\notin\ZZ_{\geq -1}$ (see \leref{psi1}).
   The restriction $|\beta_3-\beta_1| < |\beta_3-\beta_2|$ and the stronger one
    $|\beta_3-\beta_1| = |\beta_3-\beta_2|$ but $\R (\alpha_2-\alpha_1) > -1$ are 
    necessary and sufficient conditions for absolutely convergence of the number
    series
    $$\,
       F(\alpha; \beta; n)= \sum_{p=0}^{\infty}
         \frac{(2+n+\alpha_3-\alpha_2)^{(p)}}{(4+n+\alpha_3-\alpha_1)^{(p)}}
         \left(\frac{\beta_3-\beta_1}{\beta_3-\beta_2}\right)^p.
   \,$$
 This convergence ensures summability of the power series $\hat{\psi}(x)$.
 The investigation of the Stokes phenomenon of equation \eqref{initial} under relaxed
 conditions on the absolutely convergence of the series $F(\alpha; \beta; n)$ is a direction of our forther
 researches. We also postpone for further researches the general study on the
 relation between the iterated integrals and 1-sums of the products of two and
 more 1-summable power series. 

     The article is organized as follows: In the next section we briefly introduce the
     notion of Stokes matrices, as well as, the basics of
     summability theory and its application to ordinary differential equations. We also
     introduce a global funamental matrix with respect to which we are going to compute the Stokes matrices.
      In Section 3 we build a formal
     fundamental matrix at the origin. Then in Section 4 we lift this formal matrix solution to an actual fundamental
     matrix at the origin. In Section 5 we present the explicit computation of the Stokes
     matrices relevant to the formal matrix solution built in Section 3. 			
   In the last section we derive   the 1-sum
   of the product $\hat{\omega}(x)$.

	   
		\section{ Preliminaries }

		\subsection{The irregular singularity and summability theory}

    In this paragraph we briefly recall some definitions and facts from the summability theory
    (in the sense of Borel-Laplace summation)
    needed to build an actual fundamental matrix at the origin of equation \eqref{initial}.
    We follow the works of Balser  \cite{WB}, Loday-Richaud \cite{L-R}, Ramis \cite{R1, R2},
    Sauzin \cite{Sa, Sa-1}, restricting ourselves to the 1-summable case.
     
		  We denote by $\CC[[x]]$ the field of formal power series, i.e.
		$$\,
		  \CC[[x]]=\left\{ \sum_{n=0}^{\infty} f_n\,x^n\,|\,f_n\in\CC, \,n\in\NN\right\}.
		\,$$
		Equipping $\CC[[x]]$ with the natural derivation $\partial=\frac{d}{d x}$
		such that
		$$\,
		  \partial (\psi\,\varphi)=\partial(\psi)\,\varphi + \psi\,\partial(\varphi),\quad
		  \psi, \varphi\in\CC[[x]],
		  \,$$
		  we make $\CC[[x]]$ a differential algebra.
        All angular directions and sectors are defined on the Riemann surface of the
        natural logarithm.

     \bde{sec}
      1.\,An open sector $S$ is a set of the form
        $$\,
         S=S(\theta, \alpha, \rho)=\{\,x=r\,e^{i \delta}\,|\,
         0 < r < \rho,\,\theta-\alpha/1 < \delta < \theta + \alpha/2\,\},
        \,$$ 
        where $\theta$ is an arbitrary real number (the bisector of $S$), $\alpha$ is a
        positive real (the opening of $S$) and $\rho$ is either a positive number or
        $+\infty$ (the radius of $S$).\\
      2.\,A closed sector $\bar{S}$ is a set of the form
       $$\,
         \bar{S}=\bar{S}(\theta, \alpha, \rho)=\{\,x=r\,e^{i \delta}\,|\,
         0 < r \leq \rho,\,\theta -\alpha/2 \leq \delta \leq \theta+\alpha/2\,\}
       \,$$  
       with $\theta$ and $\alpha$ as before, but where $\rho$ is a positive real number
       (never equal to $+\infty$).
     \ede

    \bde{asy}
		  Let the function $f(x)$ be holomorphic on an open sector $S=S(\theta, \alpha, \rho)$.
		Let $\hat{f}(x)=\sum_{n=0}^{\infty} f_n\,x^n \in\CC[[x]]$. 
		We will say that $f(x)$ is asymptotic to $\hat{f}(x)$ in Gevrey order 1 sense on $S$,
		if  for every closed subsector $\bar{W} \subset S$  there exist  constants $C_ W > 0$
		and $A_W > 0$  such that the following estimate holds for all $N\in\NN$ and $x\in W$
			 $$\,
			   \left| f(x) - \sum_{n=0}^{N-1} f_n\,x^n \right| \leq 
			   C_W\,A_W^N\,N!\,|x|^N.
			\,$$
		\ede

       Among formal power series $\CC[[x]]$ we distinguish the formal power series
       of Gevrey order 1.
       
       \bde{gevrey}
         A formal power series $\hat{f}(x)=\sum_{n=0}^{\infty} f_n\,x^n$ is said to
         be of Gevrey order 1 if there exist two positive constants $C, A > 0$ such that
          $$\,
             |f_n| \leq C\,A^n\,n! \quad
             \textrm{for every}\quad n\in\NN.
          \,$$
       \ede 
      The constants $A$ and $C$ do not depend on $n$.
      We denote by $\CC[[x]]_1$ the set of all formal power series of Gevrey order 1.
      This set is a differential sub-algebra of $\CC[[x]]$ stable under product, derivation
      and composition \cite{L-R}.

   \bde{bo1}
     The formal Borel transform $\hat{\B}_1$ of order 1 of a formal power seies
		$\hat{f}(x)=\sum_{n=0}^{\infty} f_n\,x^n$ is called the formal series
		  $$\,
			  (\hat{\B}_1\,\hat{f}) (\zeta)= \sum_{n=0}^{\infty}
				\frac{f_n}{n!}\,\zeta^n.
			\,$$
    \ede
		
       \ble{gevrey}(\cite{R2})
         Let $\hat{f}(x)\in\CC[[x]]$ and let $f(\zeta)=(\hat{\B}_1 \hat{f})(\zeta)\in\CC[[\zeta]]$.
         Then $f(\zeta)$ converges in a neighborhood of the origin $\zeta=0$
         if and only if $\hat{f}(x)\in\CC[[x]]_1$.
       \ele	 

   Together with the formal Borel transform we consider the Laplace transform.
    \bde{laplace}
      Let $f(\zeta)$ be analytic and of exponential size at most 1 at $\infty$,
    i.e. $|f(\zeta)| \leq A \exp(B|\zeta|), \zeta\in\theta$ along any direction
    $\theta$ from 0 to $+\infty e^{i \theta}$. Then the integral
     $$\,
       (\L_{\theta} f)(x)=\int_0^{+ \infty e^{i \theta}}
       f(\zeta)\exp\left(-\frac{\zeta}{x}\right) d \left(\frac{\zeta}{x}\right)
     \,$$
     is said to be the Laplace complex transform $\L_{\theta}$ of order 1 in the
     direction $\theta$ of $f$.
    \ede
    
 Now we can give the  definition of the 1-summable power series in a direction.
  
    \bde{summable}
      The formal power series $\hat{f}(x)=\sum_{n=0}^{\infty} f_n x^n$ is 1-summable
      (or Borel summable) in the direction $\theta$ if there exist an open sector $V$
      bisected by $\theta$ whose opening is $> \pi$ and a holomorphic function $f(x)$
      on $V$ such that for every non-negative integer $N$
       $$\,
        \left| f(x) - \sum_{n=0}^{N-1} f_n x^n\right| \leq
        C_{V_1}\,A^N_{V_1}\,N!\,|x|^N
       \,$$
       on every closed subsector $\bar{V}_1$ of $V$ with constants $C_{V_1}, A_{V_1} > 0$
       depending only on $V_1$. The function $f(x)$ is called the 1-sum (or Borel sum)
       of $\hat{f}(x)$ in the direction $\theta$.
    \ede
        
        The formal power series $\hat{f}(x)$ is said to be 1-summable, if it is 1-summable
        in all but a finite number of directions.
        
     It turns out that we can verify if a Gevrey series of order 1 is  1-summable
     using the Borel and Laplace transforms. 
     
      \bpr{ver}(\cite{dP-S})
      Let $\hat{f}(x)\in\CC[[x]]_1$ and let $\theta$ be a direction. The following
      are equivalent:
        \begin{enumerate}
        	\item\,
        	$\hat{f}(x)$ is 1-summable in the direction $\theta$.
        	\item\,
        The convergent power series $(\hat{\B}_1 \hat{f})(\zeta)$ has an
        analytic continuation $h$ in a full sector
        $\{\zeta \in\CC\,|\, 0 < |\zeta| < \infty,\,|arg(\zeta) - \theta| < \epsilon\}$.
        In addition, this analytic continuation has exponential growth of order $\leq 1$
        at $\infty$ on this sector, i.e. $|h(\zeta)| \leq A\,\exp(B|\zeta|)$. In this case
        $f=\L_{\theta} (h)$ is its 1-sum. 	
        \end{enumerate}
        \epr
   Note that if the series $\hat{f}(x)\in\CC[[x]]_1$ is 1-summable in the direction $\theta$
   and $f_{\theta}(x)$ is its 1-sum in this direction then there exist an open sector $V$
   bisected by $\theta$ whose openinig is $> \pi$ such that $f_{\theta}(x)$ is asymptotic
   on $V$  to the series $\hat{f}(x)$ in Gevrey order 1 sense.
   Note also that if $\hat{f}(x)$ is convergent then its 1-sum coincides with $\hat{f}(x)$.

  We finish this paragraph introducing the differential algebra of 1-summable power series.

  Following Balser \cite{WB} we denote by $S_{1, \theta}(\hat{f})(x)$ the 1-sum of the series $\hat{f}(x)$ in the
   direction $\theta$. Denote by $\CC \{x\}_1$ the set of all 1-summable series $\hat{f}(x)$
   in all directions $\theta$ but in the same finitely many directions $\theta_1, \ldots, \theta_n$.
   This set is a differential algebra over $\CC$. In particular Theorem 2 in paragraph 3.3 in
   \cite{WB} states,

    \bth{algebra}(\cite{WB})
     For fixed, but arbitrary non-singular direction $\theta$ we have:
      \begin{enumerate}
      	\item\,
      	If $A, B\in\CC$ and $\hat{f}(x), \hat{g}(x)\in\CC\{x\}_1$ then
      	$A\,\hat{f}(x) + B\,\hat{g}(x),\,\hat{f}(x)\,\hat{g}(x)\in\CC\{x\}_1$, and
      	\ben
      	  S_{1, \theta}(A\,\hat{f} + B\,\hat{g})(x)  &=&
      	  A\,S_{1, \theta}(\hat{f})(x) + B\,S_{1, \theta}(\hat{g})(x),\\[0.15ex]
      	  S_{1, \theta}(\hat{f}\,\hat{g})(x)  &=&
      	  (S_{1, \theta}(\hat{f})(x)\,S_{1, \theta}(\hat{g})(x)).
      	 \een
      	 \item\,
      	 If $\hat{f}(x)\in\CC\{x\}_1$ then $\hat{f}'(x), \,\int_0^x \hat{f}(w)\,d w\in\CC\{x\}_1$,
      	 and
      	 \ben
      	  S_{1, \theta}(\hat{f}')(x)  &=&
      	  \frac{d}{d x} S_{1, \theta}(\hat{f})(x),\\[0.19ex]
      	  S_{1, \theta}\left(\int_0^x \hat{f}(w)\, d w\right)  &=&
      	  \int_0^x S_{1, \theta}(\hat{f})(w)\, d w.
      	 \een
      	 \item\,
      	 If $\hat{f}(x)\in\CC\{x\}_1$ has non-zero constant term, then
      	 $1/\hat{f}(x)\in\CC\{x\}_1$, and
      	 $$\,
      	   S_{1, \theta}(1/\hat{f})(x)=1/(S_{1, \theta}(\hat{f})(x)),
      	 \,$$
      	 wherever the right hand side is defined.
      	 \item\,
      	 If $\hat{f}(x)\in\CC\{x\}_1$ and $p\in\NN$, then $\hat{f}(x^p)\in\CC\{x\}_{p, \theta/p}$,
      	 and
      	 $$\,
      	   S_{p, \theta/p}\left(\hat{x}(x^p)\right)=S_{1, \theta}(\hat{f})(x^p).
      	   \,$$
      	   \item\,
      	 If $\hat{f}(x)\in\CC\{x\}_1$ has zero constant term, then
      	 $x^{-1}\,\hat{f}(x)\in\CC\{x\}_1$, and
      	 $$\,
      	   S_{1, \theta}\left(x^{-1}\,\hat{f}(x)\right)=x^{-1}\,S_{1, \theta}(\hat{f})(x).
      	 \,$$  
      \end{enumerate}    
    \ethe

      \subsection{Global, formal, actual fundamental matrices and Stokes matrices}
      
      In \cite{St} we proved that the equation \eqref{initial} possesses a 
      global fundamental matrix whose off-diagonal elements are expressed in terms
      of iterated integrals. In the present paper we are going to use
      the same fundamental matrix with respect to which we compute the Stokes matrices. 
      Theorem 2.3 in \cite{St} states
      
      \bth{global}
      Assume that $\beta_j$'s are distinct. Then the equation \eqref{initial}
      admits a global fundamental matrix $\Phi(x)$ in the form
      \be\label{global-m}
      \Phi(x)=\left(\begin{array}{ccc}
      	\Phi_1(x)   &\Phi_{12}(x)   &\Phi_{13}(x)\\[0.1ex]
      	0              &\Phi_2(x)      &\Phi_{23}(x)\\[0.1ex]
      	0              &0                     &\Phi_3(x)
      \end{array}
      \right),    
      \ee
      where the diagonal elements $\Phi_j(x), j=1, 2, 3$ are the solutions of
      the equations $L_j u=0$. The off-diagonal elements are defined as
      follows,
      \be\label{global-el}
      \Phi_{12}(x) &=&
      \Phi_1(x)\int_{\gamma_1(x)}
      \frac{\Phi_2(t)}{\Phi_1(t)} \d t,\\[0.3ex]
      \Phi_{23}(x) &=&
      \Phi_2(x)\int_{\gamma_2(x)}
      \frac{\Phi_3(t)}{\Phi_2(t)} \d t,\nonumber\\[0.3ex]
      \Phi_{13}(x) &=&
      \Phi_1(x)\int_{\gamma_3(x)}
      \frac{\Phi_2(t)}{\Phi_1(t)}
      \left( \int_{\gamma_2(t)}
      \frac{\Phi_3(t_1)}{\Phi_2(t_1)} \d t_1\right)\d t. \nonumber     
      \ee
      The paths $\gamma_j(x)$  are taken  in such a way that
      the matrix $\Phi(x)$ is a fundamental matrix
      solution of the  equation \eqref{initial}. In particular,
      the path $\gamma_1(x)$ is a path from 0 to $x$ approaching 0 in  the direction $\arg (\beta_2-\beta_1)$,
      the path $\gamma_2(x)$ is a path from 0 to $x$ approaching 0  in the direction $\arg (\beta_3-\beta_2)$
      and the path $\gamma_3(x)$ is a path from 0 to $x$ approaching 0  in the direction $\arg (\beta_3-\beta_1)$.
      \ethe

    This choice of the paths of integration $\gamma_j(x)$ implies that the
    irregular point $x=0$ is a non-resonant irregular singularity.
      
      As in \cite{St} a dierct corollary of \thref{global} gives a global fundamental
      set of solutions of the equation \eqref{initial}.
      
      \bpr{global-fss}
      Assume that $\beta_j$'s are distinct. Then the
      equation  \eqref{initial} possesses a global fundamental set of solutions of the form
      \ben
      \Phi_1(x),\quad    \Phi_1(x)\int_{\gamma_1(x)}
      \frac{\Phi_2(t)}{\Phi_1(t)} \d t,\quad
      \Phi_1(x)\int_{\gamma_3(x)}
      \frac{\Phi_2(t)}{\Phi_1(t)}
      \left( \int_{\gamma_2(t)}
      \frac{\Phi_3(t_1)}{\Phi_2(t_1)} \d t_1\right)\d t.
      \een
      
      \epr

      On the other hand as  an equation with  an irregular singularity at
      the origin of Poincar\'{e} rank 1, the equation \eqref{initial} admits
      an unique formal fundamental matrix $\hat{\Phi}(x)$ at the origin
      in the form of the theorem of Hukuhara-Turrittin \cite{W} 
      \be\label{H-T}
      \hat{\Phi}(x)=\hat{H}(x)\,x^{\Lambda}\,\exp\left(-\frac{Q}{x}\right).
      \ee
      For a non-resonant singularity $x=0$  
       the matrices $\Lambda$ and $Q$ are diagonal ones and
      \be\label{l-q}
      \Lambda=\diag(\alpha_1, \alpha_2, \alpha_3),\quad
      Q=\diag(\beta_1, \beta_2, \beta_3).
      \ee       
      Here the matrix-functions $x^{\Lambda}$ and $\exp(-Q/x)$ must be regarded
      as formal functions. The elements of the matrix $\hat{H}(x)$ are formal
      power series. Once having a formal fundamental matrix we can introduce the
      so called formal monodromy $\hat{M}$.
      
      \bde{for-mon}
      The formal monodromy matrix $\hat{M}$ relative to the formal solution
      \eqref{H-T} is defined by
      $$\,
      \hat{\Phi}(x.e^{2 \pi\,i})=\hat{\Phi}(x)\,\hat{M}.
      \,$$
      In particular,
      \ben
      \hat{M}=e^{2 \pi\,i\,\Lambda}=\left(\begin{array}{ccc}
      	e^{2 \pi\,i\,\alpha_1}  &0                       &0\\
      	0                   	&e^{2 \pi\,i\,\alpha_2}  &0   \\
      	0	                    &0                       &e^{2 \pi\,i\,\alpha_3} 
      \end{array}
      \right).    
      \een
      \ede
      
      The summability theory allows us to relate to the formal fundamental matrix
      $\hat{\Phi}(x)$ an actual fundamental matrix $\Phi(x)$ at the origin.
      More precisely
      
      \bth{actual}(Hukuhara-Turrittin-Martinet-Ramis)(\cite{JM-JR, M})
      The entries of the matrix $\hat{H}(x)$ in \eqref{H-T} are 1-summable
      in every non-singular direction $\theta$. If we denote by $H_{\theta}(x)$ the
      1-sum of $\hat{H}(x)$ along $\theta$ obtained from $\hat{H}(x)$ by a
      Borel-Laplace transform, then $\Phi_{\theta}(x)=H_{\theta}(x)\,x^{\Lambda}\,e^{-Q/x}$
      is an actual fundamental matrix at the origin of the equation \eqref{initial}.
      \ethe

     The non-zero elements $\Phi_{ij}(x),\,i \leq j$ of the formal fundamental matrix have 
     the form $c_j\,x^{\alpha_j+A_j}\,e^{-\frac{\beta_j}{x}}\,\hat{h}_{ij}(x)$,
     where $c_j\in\CC, A_j\in\NN_0$ and $\hat{h}_{ij}(x)$ are either convergent
     near the origin or diveregent power series. Classicaly we define to each divergent
     series $\hat{h}_{ij}(x)$ singular directions $\theta,\,0 \leq \theta < 2 \pi$, of the equation \eqref{initial}
     as the bisectors of any maximal angular sector where $\R (\frac{\beta_k-\beta_j}{x}) < 0$
     for some $k=1, 2, 3$. The results of the previous \cite{St} and the present work
     allows us to redefine as the unique admisible singular directions of an reducible
     third-order equation \eqref{initial} the directions
     $\theta_1=\arg(\beta_1-\beta_2),\,\theta_2=\arg(\beta_1-\beta_3)$ and $\theta_3=\arg(\beta_2-\beta_3)$.
      
       Let $\theta$ be a singular direction for the equation \eqref{initial}.
       Let $\theta^+=\theta+\epsilon$ and $\theta^-=\theta-\epsilon$, where $\epsilon > 0$
       is a small number, be two non-singular neighboring directions of the
       singular direction $\theta$. Denote by $\Phi^+_{\theta}$ and $\Phi^-_{\theta}$
       the actual fundamental matrices of the equation \eqref{initial} defined in the sense
       of \thref{actual}. Then
       
       \bde{stokes}
       With respect to a given formal fundamental matrix $\hat{\Phi}(x)$ the Stokes matrix 
       $St_{\theta}\in GL_3(\CC)$ corresponding to the singular direction $\theta$
       is defined by
       $$\,
       St_{\theta}=(\Phi^+_{\theta}(x))^{-1}\,\Phi^-_{\theta}(x),\quad
       \hat{M}\,St_{\theta}\,\hat{M}^{-1}=St_{\theta + 2 \pi}.
       \,$$
       
       \ede

  \subsection {The infinity point}
  
   The point $x=\infty$ is a regular point for the equation \eqref{initial}. Its
   characteristic exponents $\rho^{\infty}_i,\,i=1, 2, 3$ depend only on
   the parameters $\alpha_i,\,i=1, 2, 3$. In particular, they are
   $$\,
     \rho^{\infty}_1=-\alpha_1,\quad \rho^{\infty}_2=-\alpha_2-1,\quad
     \rho^{\infty}_3=-\alpha_3-2.
   \,$$
   In \cite{St} we showed that when $\alpha_1=0, \alpha_2=-2, \alpha_3=-4,
   \beta_1=1, \beta_2=2, \beta_3=0$ the point $x=\infty$ is an ordinary point.
   There are more values of the parameters for which the point $x=\infty$ becomes
   an ordinary point. For example when $\alpha_1=0, \alpha_2=\alpha_3=-3,
   \beta_3=-\beta_1=\beta_2$ the infinity point is an ordinary point too.
   In this case the path $\gamma_2(x)$ can be taken in the real positive
   axis $\RR_{+}$.

	The monodromy matrix $M_0=M^{-1}_{\infty}$ at $x=0$ (with respect to a given
	fundamental matrix) is usually called the topological or the actual monodromy
	of the equation \eqref{initial}. The topological monodromy $M_0$
	is conjugated to the product $\hat{M} St_{\theta_1} St_{\theta_2} St_{\theta_3}\in GL_3(\CC)$
	where $0 \leq \theta_1 < \theta_2 < \theta_3 < 2 \pi$ are the singular directions
	of the equation \eqref{initial}, \cite{dP-S}.

   \section{Formal matrix solution at the origin }

  For $n=0, 1, 2, \ldots$ we denote by the symbol $(a)^{(n)}$ the rising factorial
	$$\,
	   (a)^{(n)}=a\,(a+1)\,(a+2)\ldots (a+n-1),\quad
		(a)^{(0)}=1.
	\,$$
  We also denote by $\NN_0$ the set $\NN \cup \{0\}$ and by $\ZZ_{\geq m}$ (resp. $Z_{\leq m}$)
	the set of all integers $\geq m$ (resp. $\leq m$). Throughout this paper
	we denote by $\Gamma(z)$ the Euler's Gamma function.

	In this section we build a formal matrix solution at the origin of equation 
	\eqref{initial}. We start with the construction of the formal expressions at the origin
	   of the elements $\Phi_{12}(x)$ and $\Phi_{23}(x)$
	   of the formal fundamental matrix.
	   These two functions appear as a particular solution at the origin
	   of the following 
		 family of non-homogeneous first-order linear equations 	
	        \be\label{e1}
	        y'(x)=\left(\frac{\alpha_i}{x} + \frac{\beta_i}{x^2}\right)\,y(x) +
	        x^{\alpha_j}\,e^{-\frac{\beta_j}{x}}.
	        \ee	
	  The next lemma 
	 describes a particular solution at  $x=0$ of equation \eqref{e1}.
	
	\ble{l1}
	  Assume that $\beta_j$'s are distinct and $\alpha_i, \alpha_j\in\CC$. Then
		the equation \eqref{e1}
			admits a particular solution at $x=0$ in the form
			\be\label{f1}
			  y(x)=x^{\alpha_j}\,e^{-\frac{\beta_j}{x}}\,\hat{\varphi}_{ij}(x).
			\ee
			The function $\hat{\varphi}_{ij}(x)$ is defined as follows:
			\begin{enumerate}
			  \item\,
				If $\alpha_j-\alpha_i \neq -2, -3, -4, \ldots $, then
				\be\label{s1}
				  \frac{x^2}{\beta_j-\beta_i}\, \hat{\varphi}_{ij}(x)  &=&
					\frac{x^2}{\beta_j-\beta_i}
					\sum_{n=0}^{\infty} (-1)^n
					\frac{(2+\alpha_j-\alpha_i)^{(n)}}{(\beta_j-\beta_i)^n}\,x^n\\[0.3ex]
					                        &=&
					\frac{x^2}{\beta_j-\beta_i} \sum_{n=0}^{\infty} b_n\,x^n\nonumber
				\ee	
					is a divergent power series.
					\item\,
					If $\alpha_j-\alpha_i=-\alpha$ for some $\alpha=2, 3, 4, \ldots$, then
					\be\label{s2}
					  \frac{x^2}{\beta_j-\beta_i}\, \hat{\varphi}_{ij}(x)
						   &=&
					\frac{x^2}{\beta_j-\beta_i}\Big[
							1 + \frac{\alpha-2}{\beta_j-\beta_i}\,x + \frac{(\alpha-3)(\alpha-2)}{(\beta_j-\beta_i)^2}\,x^2
							+ \cdots\\[0.4ex]
						  &+&
							\frac{1.2.3\ldots(\alpha-3)(\alpha-2)}{(\beta_j-\beta_i)^{\alpha-2}}\,x^{\alpha-2}\Big]\nonumber	
					\ee
				is a finite and therefore analytic near the origin function.
			\end{enumerate}
	\ele
	
	\proof
	The proof is straightforward.\qed

In the next step we build a formal representation of $\Phi_{13}(x)$ from \eqref{global-el}  at the origin.
 In particular,
starting from the iterated integral
$$\,
\Phi_{13}(x)=x^{\alpha_1} e^{-\frac{\beta_1}{x}}
\int_0^x t^{\alpha_2-\alpha_1} e^{-\frac{\beta_2-\beta_1}{t}}
\left(\int_0^t t_1^{\alpha_3-\alpha_2} e^{-\frac{\beta_3-\beta_2}{t_1}}\,d t_1
\right) d t,
\,$$
where the inner integral is taken in the direction $\arg (\beta_3-\beta_2)$,
and the outer is taken in the direction $\arg (\beta_3-\beta_1)$, in this and the next section
we will transform $\Phi_{13}(x)$ into an appropriate function to compute
the Stokes multipliers. 

From \leref{l1} it follows that when  $\beta_3 \neq \beta_2$ then
$$\,
\int_0^x t^{\alpha_3-\alpha_2} e^{-\frac{\beta_3-\beta_2}{t}} d t=
\frac{x^{\alpha_3-\alpha_2+2} e^{-\frac{\beta_3-\beta_2}{x}}}{\beta_3-\beta_2}
\sum_{s=0}^{\sigma} (-1)^s
\frac{(2+\alpha_3-\alpha_2)^{(s)}}{(\beta_3-\beta_2)^s} x^s,
\,$$			
where the integral is taken in the direction $\arg (\beta_3-\beta_2)$.
Putting this formula in the place of the inner  integral the function $\Phi_{13}(x)$
is transformed into
\ben
\Phi_{13}(x) &=&
\frac{x^{\alpha_1} e^{-\frac{\beta_1}{x}}}{\beta_3-\beta_2}
\int_0^x t^{\alpha_3-\alpha_1+2} e^{-\frac{\beta_3-\beta_1}{t}}
\left(\sum_{s=0}^{\sigma} (-1)^s \frac{(2+\alpha_3-\alpha_2)^{(s)}}{(\beta_3-\beta_2)^s} t^s\right)
d t\\[0.25ex]
&=&
\frac{x^{\alpha_1} e^{-\frac{\beta_1}{x}}}{\beta_3-\beta_2}\,I.  
\een			
Here $\sigma=\alpha_2-\alpha_3-2$ if $\alpha_3-\alpha_2\in\ZZ_{\leq -2}$, and $\sigma=\infty$
otherwise.

The next lemma  gives a formal representation of the function $\Phi_{13}(x)$ at the origin.

\ble{l3}
Assume that $\beta_j$'s are distinct. 
Then the function $\Phi_{13}(x)$ is represented at the origin as follows,
\begin{enumerate}
	\item\,
	If $\alpha_3-\alpha_1\notin\ZZ_{\leq -4}$ then
	\ben
	\Phi_{13}(x)  =
	\frac{x^{\alpha_3+4}\,e^{-\frac{\beta_3}{x}}}{(\beta_3-\beta_2)\,(\beta_3-\beta_1)}
	\left[F(\alpha; \beta; 0)\sum_{s=0}^{\infty}
	(-1)^s \frac{(4+\alpha_3-\alpha_1)^{(s)}}{(\beta_3-\beta_1)^s} x^s
	- \sum_{n=0}^{\sigma} a_n\,x^n\right],
	\een
	where the  number series $F(\alpha; \beta; 0)$ is given by 
	\be\label{F0}
	F(\alpha; \beta; 0)=\sum_{p=0}^{\sigma}
	\frac{(2+\alpha_3-\alpha_2)^{(p)}}{(4+\alpha_3-\alpha_1)^{(p)}}
	\left(\frac{\beta_3-\beta_1}{\beta_3-\beta_2}\right)^p.
	\ee
	The coefficients $a_n$ are defined by
	\be\label{an}
	a_n=b_n\,F(\alpha; \beta; n) - b_n,
	\ee
	where
	$$\,
	F(\alpha; \beta; n)=
	\sum_{p=0}^{\sigma}
	\frac{(2+n+\alpha_3-\alpha_2)^{(p)}}{(4+n+\alpha_3-\alpha_1)^{(p)}}
	\left(\frac{\beta_3-\beta_1}{\beta_3-\beta_2}\right)^p
	\,$$
	and 
	$$\,
	b_n=(-1)^n\,\frac{(2 + \alpha_3-\alpha_2)^{(n)}}{(\beta_3-\beta_2)^n}.
	\,$$
	Here  $\sigma=\alpha_2-\alpha_3-2$ if $\alpha_3-\alpha_2\in\ZZ_{\leq -2}$,
	and $\sigma=\infty$ otherwise.
	
	\item\,
	If $\alpha_3-\alpha_1=-\alpha$ for $\alpha\in\ZZ_{\geq 4}$ but
	$\alpha_3-\alpha_2$ and $\alpha_2-\alpha_1$ do not belong to $\ZZ_{\leq -2}$
	simultaneously then
	\ben
	\Phi_{13}(x)  &=&
	\frac{x^{\alpha_3+4}\,e^{-\frac{\beta_3}{x}}}{(\beta_3-\beta_2)\,(\beta_3-\beta_1)}
	\sum_{l=0}^{\alpha-4} b_l\,x^l\left(\sum_{s=0}^{\alpha-l-4}
	(-1)^s \frac{(4+l+\alpha_3-\alpha_1)^{(s)}}{(\beta_3-\beta_1)^s}\,x^s\right)\\[0.4ex]
	&+&
	\frac{b_{\alpha-3}\,x^{\alpha_1+1}\,e^{-\frac{\beta_3}{x}}}{(\beta_3-\beta_2)\,(\beta_3-\beta_1)}
	\left[\tilde{F}(\alpha; \beta; 0)
	\sum_{s=0}^{\infty}
	(-1)^s \frac{s!}{(\beta_3-\beta_1)^s} x^s
	- \sum_{n=0}^{\sigma} a_n\,x^n\right],
	\een
	where 
	$$\,
	\tilde{F}(\alpha; \beta; 0)=\sum_{p=0}^{\infty}
	\frac{(\alpha_1-\alpha_2-1)^{(p)}}{(1)^{(p)}}
	\left(\frac{\beta_3-\beta_1}{\beta_3-\beta_2}\right)^p=
	\left(\frac{\beta_1-\beta_2}{\beta_3-\beta_2}\right)^{\alpha_2-\alpha_1+1}.
	\,$$
	The coefficients $a_n$ are defined by
       \be\label{ant}
	a_n=c_n\,\tilde{F}(\alpha; \beta; n) - c_n,
	     \ee
	where
	$$\,
	\tilde{F}(\alpha; \beta; n)=\sum_{p=0}^{\infty}
	\frac{(\alpha_1-\alpha_2-1+n)^{(p)}}{(1+n)^{(p)}}
	\left(\frac{\beta_3-\beta_1}{\beta_3-\beta_2}\right)^p
	\,$$
	and
	$$\,
	c_n=(-1)^n \frac{(\alpha_1-\alpha_2-1)^{(n)}}{(\beta_3-\beta_2)^n}.
	\,$$
	Here $\sigma=\alpha_2-\alpha_1$ when $\alpha_2-\alpha_1\in\ZZ_{\geq -1}$ and
	$\sigma=\infty$ otherwise. When  $\alpha_2-\alpha_1=-1$ the sum $\sum_{n=0}^{-1} a_n\,x^n$ is equal to zero.

	\item\,
	If $\alpha_3-\alpha_2,\,\alpha_2-\alpha_1\in\ZZ_{\leq -2}$ then
	$$\,
	\Phi_{13}(x)=
		\frac{x^{\alpha_3+4}\,e^{-\frac{\beta_3}{x}}}{(\beta_3-\beta_2)\,(\beta_3-\beta_1)}
		\sum_{l=0}^{\alpha_2-\alpha_3-2} b_l\,x^l
		\left(\sum_{s=0}^{\alpha_1-\alpha_3-4-l}
		(-1)^s \frac{(4+l+\alpha_3-\alpha_1)^{(s)}}{(\beta_3-\beta_1)^s}\,x^s\right).
	\,$$

\end{enumerate}
\ele

\proof

The integral $I$ can be written as
$$\,
I=\sum_{n=0}^{\sigma} I_n,\quad
I_n=b_n\int_0^x t^{\alpha_3-\alpha_1+2+n}\,e^{-\frac{\beta_3-\beta_1}{t}}\,d t,
\,$$
where
$$\,
b_n=(-1)^n\,\frac{(2 + \alpha_3-\alpha_2)^{(n)}}{(\beta_3-\beta_2)^n}.
\,$$
Here $\sigma=\alpha_2-\alpha_3-2$ if $\alpha_3-\alpha_2\in\ZZ_{\leq -2}$, and
$\sigma=\infty$ otherwise.
For $\alpha_3-\alpha_1\notin\ZZ_{\leq -4}$ every integral $I_n$ is expressed
as an infinite power series
\ben
I_n &=&
\frac{b_n\,x^{\alpha_3-\alpha_1+4+n}\,e^{-\frac{\beta_3-\beta_1}{x}}}{\beta_3-\beta_1}
\sum_{s=0}^{\infty} (-1)^s
\frac{(4+n+\alpha_3-\alpha_1)^{(s)}}{(\beta_3-\beta_1)^s}\,x^s\\[0.4ex]
&=&
\frac{(-1)^n\,b_n\,x^{\alpha_3-\alpha_1+4} e^{-\frac{\beta_3-\beta_1}{x}}} {\beta_3-\beta_1}
\frac{(\beta_3-\beta_1)^n}{(4+\alpha_3-\alpha_1)^{(n)}}
\left[S_0(x)-\sum_{s=0}^{n-1} (-1)^s
\frac{(4+\alpha_3-\alpha_1)^{(s)}}{(\beta_3-\beta_1)^s}\,x^s\right],     
\een
where 
$$\,
S_0(x)=\sum_{s=0}^{\infty} (-1)^s \frac{(4+\alpha_3-\alpha_1)^{(s)}}{(\beta_3-\beta_1)^s} x^s.
\,$$
Then $I$ becomes
$$\,
I=\frac{x^{\alpha_3-\alpha_1+4}\,e^{-\frac{\beta_3-\beta_1}{x}}}{\beta_3-\beta_1}
\left[S_0(x)\,F(\alpha; \beta; 0) - \sum_{n=0}^{\sigma} a_n\,x^n\right],
\,$$
where $F(\alpha; \beta; 0)$ and $a_n$ are defined by \eqref{F0} and \eqref{an},
respectively.

Let $\alpha_3-\alpha_1=-\alpha$ for $\alpha\in\ZZ_{\geq -4}$. Then every integral
$$\,
I_l=b_l\int_0^x t^{\alpha_3-\alpha_1+2+l}\,e^{-\frac{\beta_3-\beta_1}{t}}\, d t,\quad
0 \leq l \leq -4+\alpha
\,$$
is expressed in  terms of finite series. In particular, in this case we have
$$\,
I_l=
b_l \frac{x^{\alpha_3-\alpha_1+4+l}\,e^{-\frac{\beta_3-\beta_1}{x}}}{\beta_3-\beta_1}
\sum_{s=0}^{\alpha-l-4} (-1)^s
\frac{(4+n+\alpha_3-\alpha_1)^{(s)}}{(\beta_3-\beta_1)^s}\,x^s.
\,$$
If in addition $\alpha_3-\alpha_2,\,\alpha_2-\alpha_1\in\ZZ_{\leq -2}$ then
  $$\,
    I=\sum_{l=0}^{\alpha_2-\alpha_3-2} I_l=
    \frac{x^{\alpha_3-\alpha_1+4}\,e^{-\frac{\beta_3-\beta_1}{x}}}{\beta_3-\beta_1}
    \sum_{l=0}^{\alpha_2-\alpha_3-2} b_l\,x^l\left(
    \sum_{s=0}^{\alpha-l-4} (-1)^s
    \frac{(4+n+\alpha_3-\alpha_1)^{(s)}}{(\beta_3-\beta_1)^s}\,x^s\right).
  \,$$
Let $\alpha_3-\alpha_2$ and $\alpha_2-\alpha_1$ do not belong to $\ZZ_{\leq -2}$ simultaneously.
Then for $l \geq -3+\alpha$ we have that
\ben
I_{\alpha-3+n}=(-1)^n\,b_{\alpha-3+n}
\frac{x\,e^{-\frac{\beta_3-\beta_1}{x}}}{\beta_3-\beta_1}
\frac{(\beta_3-\beta_1)^n}{(1)^{(n)}}
\left[\tilde{S}_0(x) - \sum_{s=0}^{n-1} \frac{s!}{(\beta_3-\beta_1)^s}\,x^s\right],
\een
where
$$\,
\tilde{S}_0(x)=\sum_{s=0}^{\infty} (-1)^s
\frac{s!}{(\beta_3-\beta_1)^s}\,x^s.
\,$$
Then the integral $I$ becomes
\ben
I=\sum_{l=0}^{\alpha-4} I_l
+\frac{x\,e^{-\frac{\beta_3-\beta_1}{x}}}{\beta_3-\beta_1}
\left[b_{\alpha-3}\,\tilde{F}(\alpha; \beta; 0)\,\tilde{S}_0(x) -
b_{\alpha-3}\,\sum_{n=0}^{\infty} a_n\,x^n\right],
\een
where
$$\,
\tilde{F}(\alpha; \beta; 0)=\sum_{p=0}^{\infty}
\frac{(\alpha_1-\alpha_2-1)^{(p)}}{(1)^{(p)}}
\left(\frac{\beta_3-\beta_1}{\beta_3-\beta_2}\right)^p=
\left(\frac{\beta_1-\beta_2}{\beta_3-\beta_2}\right)^{\alpha_2-\alpha_1+1}.
\,$$
The coefficients $a_n$ are defined by
$$\,
a_n=c_n\,\tilde{F}(\alpha; \beta; n) - c_n,
\,$$
where
$$\,
\tilde{F}(\alpha; \beta; n)=\sum_{p=0}^{\infty}
\frac{(\alpha_1-\alpha_2-1+n)^{(p)}}{(1+n)^{(p)}}
\left(\frac{\beta_3-\beta_1}{\beta_3-\beta_2}\right)^p
\,$$
and
$$\,
c_n=(-1)^n \frac{(\alpha_1-\alpha_2-1)^{(n)}}{(\beta_3-\beta_2)^n}.
\,$$
 In the particular case  when $\alpha_3-\alpha_2\in\ZZ_{\leq -2}$ but
$\alpha_2-\alpha_1\in\ZZ_{\geq -1}$  the integral $I$ is reduced to
   $$\,
  I=\sum_{l=0}^{\alpha-4} I_l
  +\frac{b_{\alpha-3}\,x\,e^{-\frac{\beta_3-\beta_1}{x}}}{\beta_3-\beta_1}
  \left[\tilde{F}(\alpha; \beta; 0)\,\tilde{S}_0(x) -
  \sum_{n=0}^{\alpha_2-\alpha_1} a_n\,x^n\right],
    \,$$
    since $c_n=0$ for $n \geq \alpha_2-\alpha_1+2$ and
    $\tilde{F}(\alpha; \beta; \alpha_2-\alpha_1+1) -1=0$. We also note that
    when $\alpha_2-\alpha_1=-1$ we define the series $\sum_{n=0}^{-1} a_n\,x^n=a_0=0$
    since in this case $\tilde{F}(\alpha; \beta; 0)=1$.
    
This completes the proof.
\qed

		Now we can make a formal fundamental matrix  $\hat{\Phi}(x)$   at the origin explicit.

  \bpr{fs-0}
   Assume that $\beta_j$'s are distinct. 
   Then the  equation \eqref{initial} 
   possesses an unique formal fundamental matrix $\hat{\Phi}(x)$ 
   at the origin in the form \eqref{H-T}-\eqref{l-q}
      where
      \be\label{h-0}\qquad
      \hat{H}(x)=\left(\begin{array}{ccc}
      	1   &\frac{x^2\,\hat{\varphi}_{12}(x)}{\beta_2-\beta_1}       
      	&\frac{x^4\,\hat{\psi}(x)}{(\beta_3-\beta_2) (\beta_3-\beta_1)}\\[0.35ex]
      	0   &1                           &\frac{x^2\,\hat{\varphi}_{23}(x)}{\beta_3-\beta_2}\\
      	0   &0                           &1
      \end{array}
      \right).
      \ee
   The formal function $x^2 \hat{\varphi}_{ij}(x)/(\beta_j-\beta_i)$
   is defined by the finite series \eqref{s2} if $\alpha_j-\alpha_i\in\ZZ_{\leq -2}$,
   and by the infinite formal series \eqref{s1} otherwise.
   The formal function $\hat{\psi}(x)$ 
    is defined  as follows:
     \begin{enumerate}
     \item\,
     If 
     $\alpha_3-\alpha_1\notin\ZZ_{\leq -4}$ then
           $$\,
       \hat{\psi}(x)=F(\alpha; \beta; 0)\,S_0(x) - \sum_{n=0}^{\sigma} a_n\,x^n,
      \,$$
     where $F(\alpha; \beta; 0), S_0(x)$ and $a_n$ are defined by the first part of
     \leref{l3}. Here $\sigma=\alpha_2-\alpha_3-2$ if $\alpha_3-\alpha_2\in\ZZ_{\leq -2}$,
     and $\sigma=\infty$ otherwise.
    
     \item\,
     If  $\alpha_3-\alpha_1\in\ZZ_{\leq -4}$ but 
     $\alpha_2-\alpha_1$ and $\alpha_3-\alpha_2$ do not belong to $\ZZ_{\leq -2}$ simultaneously then  
      \ben
        \hat{\psi}(x)
          &=&
        \sum_{l=0}^{\alpha_1-\alpha_3-4} b_l\,x^l \left(
        \sum_{s=0}^{\alpha_1-\alpha_3-4-l} (-1)^s\,\frac{(4+l+\alpha_3-\alpha_1)^{(s)}}
        {(\beta_3-\beta_1)^s}\,x^s\right) \\[0.4ex]
        &+& 
        x^{-(\alpha_3-\alpha_1)-3}\,b_{\alpha_1-\alpha_3-3}
        \left[\tilde{F}(\alpha; \beta; 0)\,\tilde{S}_0(x) - 
        \sum_{n=0}^{\sigma} a_n\,x^n\right],
      \een
      where
      $$\,
        b_n=(-1)^n \frac{(2+\alpha_3-\alpha_2)^{(n)}}{(\beta_3-\beta_2)^n}.
      \,$$
      The series $\tilde{F}(\alpha; \beta; 0),\,\tilde{S}_0(x)$ and the coefficients $a_n$
      are defined by the second part of \leref{l3}.
      Here $\sigma=\alpha_2-\alpha_1$ if $\alpha_2-\alpha_1\in\ZZ_{\geq 0}$,
      $\sigma=0$ if $\alpha_2-\alpha_1=-1$ and $\sigma=\infty$ otherwise.
   
   \item\,
   If $\alpha_3-\alpha_2,\,\alpha_2-\alpha_1\in\ZZ_{\leq -2}$  then     
         $$\,
         \hat{\psi}(x)=
        \sum_{l=0}^{\alpha_2-\alpha_3-2} b_l\,x^l
        \left(\sum_{s=0}^{\alpha_1-\alpha_3-4-l} (-1)^s\,
        \frac{(4+l+\alpha_3-\alpha_1)^{(s)}}{(\beta_3-\beta_1)^s}\,x^s\right).
    \,$$

  \end{enumerate}
      \epr

\section{Actual fundamental matrix at the origin}

  In this paragraph we  find the 1-sum $H(x)$ of the formal
  function $\hat{H}(x)$. To solve this problem we utilize the summability
  theory in a direction. The product of the obtained actual function
  $H(x)$ with the actual function $x^{\Lambda}\,\exp(-Q/x)$ yields
  an actual solution of the  equation \eqref{initial} at the origin, corresponding
  to the formal solution from \prref{fs-0}.

   The first lemma gives the 1-sums of
   the infinite formal power series \eqref{s1}  
   and the series  $S_0(x),\,\tilde{S}_0(x)$ inroduced by \leref{l3}.
  
  \ble{phi}
  \begin{description}
  	\item[I.]\,
  	Assume that $\beta_j$'s are distinct and that $\alpha_j-\alpha_i \neq -2, -3, -4, \ldots$.
  	Then the formal power series
  	\be\label{fps}
  	\hat{\varphi}(x)=\sum_{n=0}^{\infty}
  	(-1)^n \frac{(2+\alpha_j-\alpha_i)^{(n)}}{(\beta_j-\beta_i)^n}\,x^n
  	\ee
  	is  1-summable. For any direction $\theta \neq \arg(\beta_i-\beta_j)$ from
  	$0$ to $+\infty\,e^{i \theta}$   	the function
  	\be\label{al}
  	\varphi_{\theta}(x)=\int_0^{+\infty\,e^{i \theta}}
  	\left(1+\frac{\zeta}{\beta_j-\beta_i}\right)^{\alpha_i-\alpha_j-2}\,
  	e^{-\frac{\zeta}{x}}\,d \left(\frac{\zeta}{x}\right)
  	\ee
  	defines its 1-sum in such a direction.
  	
  	\item[II.]\,
  	Assume that $\alpha_3-\alpha_1\notin\ZZ_{\leq -4}$. Then the  formal
  	power series
  	\ben
  	\hat{\upsilon}(x)=\sum_{s=0}^{\infty}
  	(-1)^s\,
  	\frac{(4+\alpha_3-\alpha_1)^{(s)}}{(\beta_3-\beta_1)^s}\,x^s,
  	\een		
  	is  1-summable. For any direction $\theta \neq \arg (\beta_1-\beta_3)$
  	from $0$ to $+\infty\,e^{i \theta}$  
  	the function
  	\be\label{al-1}
  	\upsilon_{\theta}(x)= 
  	\int_0^{+\infty\,e^{i \theta}}
  	\left(1 + \frac{\nu}{\beta_3-\beta_1}\right)^{\alpha_1-\alpha_3-4}\,
  	e^{-\frac{\nu}{x}}\,d \left(\frac{\nu}{x}\right)
  	\ee 
  	defines its 1-sum in such a direction.
  \end{description}
  \ele

  \proof
  We focus on the proof of the first statement.
  
  Let $|2+\alpha_j-\alpha_i| \leq 1$. Then
  $$\,
  |(2+\alpha_j-\alpha_i)^{(n)}| \leq n!.
  \,$$
  Let $|2+\alpha_j-\alpha_i| > 1$. Then
  $$\,
  |(2+\alpha_j-\alpha_i)^{(n)}| \leq (|2+\alpha_j-\alpha_i|+1)^n\,n!.
  \,$$
  Therefore the formal power series $\hat{\varphi}(x)$ is of Gevrey order 1 with
  constants
  $$\,
  C=1,\quad A=\frac{1}{|\beta_j-\beta_i|}
  \,$$
  if $|2+\alpha_j-\alpha_i| \leq 1$, and
  $$\,
  C=1,\quad A=\frac{|2+\alpha_j-\alpha_i|+1}{|\beta_j-\beta_i|}
  \,$$
  if $|2+\alpha_j-\alpha_i| > 1$.
  
  As a result the formal Borel transform
  $$\,
  \varphi(\zeta)=(\hat{\B}_1\,\hat{\varphi})(\zeta)=
  \sum_{n=0}^{\infty}
  (-1)^n \frac{(2+\alpha_j-\alpha_i)^{(n)}}{(\beta_j-\beta_i)^n}\,
  \frac{\zeta^n}{n!}
  \,$$
  is an analytic function near the origin in the Borel plane.
  
  Then for any direction $\theta \neq \arg(\beta_i-\beta_j)$ from $0$
  to $+\infty\,e^{i \theta}$ the associate Laplace transform
  $$\,
  \varphi_{\theta}(x)=\int_0^{+\infty\,e^{i \theta}}
  \left(1+\frac{\zeta}{\beta_j-\beta_i}\right)^{\alpha_i-\alpha_j-2}\,
  e^{-\frac{\zeta}{x}}\,d \left(\frac{\zeta}{x}\right)
  \,$$
  defines its 1-sum in such a direction.
  
  In the same manner one can prove that 	the series 
  $$\,
  \hat{\upsilon}(x)=\sum_{s=0}^{\infty} (-1)^s \frac{(4+\alpha_3-\alpha_1)^{(s)}}{(\beta_3-\beta_1)^s}\,x^s			
  \,$$			
  is  1-summable in any direction $\theta \neq \arg (\beta_1-\beta_3)$ form 0 to
  $+\infty\,e^{i \theta}$ and the function
  $$\,
  \upsilon_{\theta}(x)=
  \int_0^{+ \infty\,e^{i \theta}}
  \left(1 + \frac{\nu}{\beta_3-\beta_1}\right)^{\alpha_1-\alpha_3-4}\,
  e^{-\frac{\nu}{x}}\,d \left(\frac{\nu}{x}\right)
  \,$$
  defines its 1-sum in such a direction. 
  
  This completes the proof.
  \qed

\bre{1-sum-1}
   When we move the direction $\theta$ continuosly the corresponding 1-sums
   $\varphi_{\theta}(x)$ (resp. $\upsilon_{\theta}(x)$) stick each other analytically
   and define a holomorphic function $\tilde{\varphi}(x)$ (resp. $\tilde{\upsilon}(x)$)
  on a sector of opening $3 \pi$. Denote by $\theta_1$ the singular direction of
      $\varphi_{\theta}(x)$ or $\upsilon_{\theta}(x)$.
    When $\pi/2 < \theta_1 \leq 3 \pi/2$ this sector is defined by
   $\theta_1-2 \pi-\pi/2 < \arg (x) <\theta_1 + \pi/2$. When $-\pi/2 < \theta_1 \leq \pi/2$
   this sector is defined by $\theta_1-\pi/2 < \arg (x) < \theta_1+2 \pi +\pi/2$
   (see paragraph 5.10 in \cite{Sa}).
  On these sectors the multivalued functions $\tilde{\varphi}(x)$ and $\tilde{\upsilon}(x)$
  define the 1-sums of the formal series $\hat{\varphi}(x)$ and $\hat{\upsilon}(x)$
  and these functions are asymptotic to these series in Gevrey order 1 sense.
  In every non-singular direction $\theta$ the multivalued functions $\tilde{\varphi}(x)$
  and $\tilde{\upsilon}(x)$ have one value $\varphi_{\theta}(x)$ and $\upsilon_{\theta}(x)$,
  respectively. Near the corresponding singular direction $\theta_1$ the functions
  $\tilde{\varphi}(x)$ and $\tilde{\upsilon}(x)$ have two different values:
  $\varphi^+_{\theta_1}(x)=\varphi_{\theta_1+\epsilon}(x),\,
  \upsilon^+_{\theta_1}(x)=\upsilon_{\theta_1+\epsilon}(x)$ and
    $\varphi^-_{\theta_1}(x)=\varphi_{\theta_1-\epsilon}(x),\,
    \upsilon^-_{\theta_1}(x)=\upsilon_{\theta_1-\epsilon}(x)$, where $\epsilon > 0$ is a
    small number.
   \qed
\ere

  The next lemma describes the behaviour of the infinite number series $F(\alpha ; \beta; n)$
  and $\tilde{F}(\alpha; \beta; n)$ for $n\in\NN_0$.
  
  \ble{S2}
  Assume that $\alpha_3-\alpha_1 \notin\ZZ_{\leq -4}$ 
  Then for every $n\in\NN_0$ the number series
  $$\,
  F(\alpha ; \beta; n)=
  \sum_{p=0}^{\infty}
  \frac{(2+n+\alpha_3-\alpha_2)^{(p)}}{(4 + n +\alpha_3-\alpha_1)^{(p)}}
  \left(\frac{\beta_3-\beta_1}{\beta_3-\beta_2}\right)^p
  $$
  and
  $$\,
  \tilde{F}(\alpha; \beta; n)=\sum_{p=0}^{\infty}
  \frac{(\alpha_1-\alpha_2-1+n)^{(p)}}{(1+n)^{(p)}}
  \left(\frac{\beta_3-\beta_1}{\beta_3-\beta_2}\right)^p
  \,$$
  satisfy the following conditions:
  \begin{enumerate}
  	\item\,
  	If $|\beta_3-\beta_1| < |\beta_3-\beta_2|$, then both $F(\alpha ; \beta; n)$
  	and $\tilde{F}(\alpha; \beta; n)$ are
  	absolutely convergent series.
  
  	\item\,
  	If $|\beta_3-\beta_1| = |\beta_3-\beta_2|$ and $\R (\alpha_1-\alpha_2) < 1$, then both
  	$F(\alpha ; \beta; n)$ and $\tilde{F}(\alpha; \beta; n)$ are absolutely convergent series.
  
  	\item\,
  	If $|\beta_3-\beta_1| = |\beta_3-\beta_2|,\,\beta_3-\beta_1 \neq \beta_3-\beta_2$ 
  	and $1 \leq \R (\alpha_1-\alpha_2) < 2$,
  	then both $F(\alpha ; \beta; n)$  and $\tilde{F}(\alpha; \beta; n)$ are conditionally  convergent series.
  
  	\item\,
  	If $|\beta_3-\beta_1| = |\beta_3-\beta_2|$ and $\R (\alpha_1-\alpha_2) \geq 2$, then
  	both $F(\alpha ; \beta; n)$ and $\tilde{F}(\alpha; \beta; n)$ are  divergent series.
  
    \item\,
    If $|\beta_3-\beta_1| > |\beta_3-\beta_2|$, then both $F(\alpha; \beta; n)$ and 
    $\tilde{F}(\alpha; \beta; n)$ are divergent series.
  \end{enumerate}
  \ele

  \proof
  The series $F(\alpha ; \beta; n)$ and $\tilde{F}(\alpha; \beta; n)$ are nothing but the hypergeometric series
  $$\,
  _2F_1(a, b ; c ; z)=\sum_{p=0}^{\infty}\frac{(a)^{(p)}\,(b)^{(p)}}{(c)^{(p)}\,p!}\,z^p
  \quad (\textrm{in short}\quad F(a, b; c; z))
  \,$$
  for $a=1, b=2+n+\alpha_3-\alpha_2$,
  $c=4+n+\alpha_3-\alpha_1$ and $z=(\beta_3-\beta_1)/(\beta_3-\beta_2)$ for the series $F$.
  For the series $\tilde{F}$ the parameters are $a=1, b=\alpha_1-\alpha_2-1+n, c=1+n$. 
  Then the statement follows from the properties of
  the hypergeometric series. In particular, if  $c$ is different from
  $0, -1, -2, -3, \ldots$, then the hypergeometric series is absolutely convergent
  for all $|z| < 1$. Note that on the assumptions of this lemma
  $c$ is different from 
  $0, -1, -2, -3, \ldots$ for both series $F$ and $\tilde{F}$.
  When $|z|=1$ we have that the hypergeometric series
  is absolutely convergent if $\R (a+b-c) < 0$;
  is conditionally convergent if $z \neq 1$ and $0 \leq \R (a+b-c) < 1$;
  is divergent if $1 \leq \R (a+b-c)$ (see 2.1.1 in \cite{HB-AE})\,.
  \qed
  
 In this paper we work with absolutely convergent number series $F$ and $\tilde{F}$
 without making other restrictions on the parameters. So, we assume that
  $|\beta_3-\beta_1| < |\beta_3-\beta_2|$. At the end of this section we show that the next results
  remain valid for these values of the parameters $\beta_j$'s for which $|\beta_3-\beta_1|=|\beta_3-\beta_2|$
  but $\beta_3-\beta_1 \neq \pm (\beta_3-\beta_2)$ provided that $\R (\alpha_2-\alpha_1) > -1$.

  Now we deal with the formal series $\sum a_n\,x^n$ 
  itroduced in  \leref{l3}.

  Consider the formal power series
  \ben
  \hat{\psi}(x)=\sum_{n=0}^{\infty} a_n\,x^n,&\quad&
  a_n=b_n\, F(\alpha; \beta; n) -
  b_n \,,\\
  \hat{\phi}(x)=\sum_{n=0}^{\infty} a_n\,x^n,&\quad&
  a_n=c_n\, \tilde{F}(\alpha; \beta; n) -
  c_n ,
  \een
  Represent $\hat{\psi}(x)$ and $\hat{\phi}(x)$ as
  $$\,
  \hat{\psi}(x)=\hat{\psi_1}(x)-\hat{\psi_2}(x),\quad
  \hat{\phi}(x)=\hat{\phi_1}(x)-\hat{\phi_2}(x),
  \,$$
  where
  \ben
  \hat{\psi_1}(x) &=&
  \sum_{n=0}^{\infty}
  \frac{(-1)^n\,(2+\alpha_3-\alpha_2)^{(n)}}{(\beta_3-\beta_2)^n}\,
  F(\alpha; \beta; n)\,x^n,\\[0.3ex]
  \hat{\phi_1}(x) &=&
  \sum_{n=0}^{\infty}
  \frac{(-1)^n\,(\alpha_1-\alpha_2-1)^{(n)}}{(\beta_3-\beta_2)^n}\,
  \tilde{F}(\alpha; \beta; n)\,x^n,\\[0.3ex]
  \hat{\psi_2}(x)  &=&
  \sum_{n=0}^{\infty}
  \frac{(-1)^n\,(2+\alpha_3-\alpha_2)^{(n)}}{(\beta_3-\beta_2)^n}\,x^n,\,\,
  \hat{\phi_2}(x)=
  \sum_{n=0}^{\infty}
  \frac{(-1)^n\,(\alpha_1-\alpha_2-1)^{(n)}}{(\beta_3-\beta_2)^n}\,x^n.
  \een

  In the following lemma we build the 1-sum of the formal series $\hat{\psi}(x)$. 
  
  \ble{psi1}
  Assume that $|\beta_3-\beta_1| < |\beta_3-\beta_2|$ provided that $\beta_j$'s
  are distinct . 
  Assume also that $\alpha_3-\alpha_1\notin\ZZ_{\leq -4}$ and
  $\alpha_3-\alpha_2\notin\ZZ_{\leq -2}$.
  Then the formal power series $\hat{\psi}(x)$ is 1-summable.
    If $\alpha_3-\alpha_2\notin\ZZ$	
  	for any directions $\theta_3 \neq \arg (\beta_2-\beta_3)$ from $0$
  	to $+\infty\,e^{i \theta_3}$ and $\theta_2 \neq \arg (\beta_1-\beta_3)$ from $0$ to
  	$+\infty\,e^{i \theta_2}$   the function
  	\ben
  	   & &
  	\psi_{\theta}(x) =
  	- \frac{\beta_3-\beta_1}{\beta_3-\beta_2}
  	\frac{\Gamma(4+\alpha_3-\alpha_1) \Gamma(\alpha_2-\alpha_3-1)}
  	{\Gamma(2+\alpha_2-\alpha_1)}
  	\left(\frac{\beta_1-\beta_2}{\beta_3-\beta_2}\right)^{\alpha_2-\alpha_3-3}\\[0.4ex]
  	 &\times&
  	\left(\frac{\beta_1-\beta_3}{\beta_1-\beta_2}\right)^{\alpha_1-\alpha_3-4}
  	\int_0^{+\infty e^{i \theta_2}}
  	\left(1+\frac{\xi}{\beta_3-\beta_1}\right)^{\alpha_1-\alpha_3-4}
  	e^{-\frac{\xi}{x}} d\left(\frac{\xi}{x}\right)\\[0.85ex]  
  	&+&
  	\frac{\beta_3-\beta_1}{\beta_1-\beta_2}
  	\sum_{s=0}^{\infty}
  	\frac{(2+\alpha_2-\alpha_1)^{(s)}}{(\alpha_2-\alpha_3-1)^{(s)}}
  	\left(\frac{\beta_3-\beta_2}{\beta_1-\beta_2}\right)^s
  	\int_0^{+\infty\,e^{i \theta_3}}
  	\left(1+\frac{\xi}{\beta_3-\beta_2}\right)^{\alpha+s}\,
  	e^{-\frac{\xi}{x}}\,d \left(\frac{\xi}{x}\right)
  	\een
  	where $\alpha=\alpha_2-\alpha_3-2$
  	defines the corresponding 1-sum in such  directions.
  
    If $\alpha_3-\alpha_2\in\ZZ_{\geq -1}$ but $\alpha_2-\alpha_1\in\ZZ_{\leq -2}$
  the 1-sum $\psi_{\theta}(x)$ has the form
       	\ben
       	   & &
       	\psi_{\theta}(x) =
       	- \frac{\beta_3-\beta_1}{\beta_3-\beta_2}
       	\frac{\Gamma(4+\alpha_3-\alpha_1) }
       	{(\alpha_2-\alpha_3-1)(\alpha_2-\alpha_3)\dots (\alpha_2-\alpha_3+m-2)}\\[0.4ex]
       	  &\times&
       	\left(\frac{\beta_1-\beta_2}{\beta_3-\beta_2}\right)^{\alpha_2-\alpha_3-3}
       	\left(\frac{\beta_1-\beta_3}{\beta_1-\beta_2}\right)^{\alpha_1-\alpha_3-4}
       	\int_0^{+\infty e^{i \theta_2}}
       	\left(1+\frac{\xi}{\beta_3-\beta_1}\right)^{\alpha_1-\alpha_3-4}
       	e^{-\frac{\xi}{x}} d\left(\frac{\xi}{x}\right)\\[0.85ex]  
       	&+&
       	\frac{\beta_3-\beta_1}{\beta_1-\beta_2}
       	\sum_{s=0}^{\alpha}
       	\frac{(2+\alpha_2-\alpha_1)^{(s)}}{(\alpha_2-\alpha_3-1)^{(s)}}
       	\left(\frac{\beta_3-\beta_2}{\beta_1-\beta_2}\right)^s
       	\int_0^{+\infty\,e^{i \theta_3}}
       	\left(1+\frac{\xi}{\beta_3-\beta_2}\right)^{\alpha+s}\,
       	e^{-\frac{\xi}{x}}\,d \left(\frac{\xi}{x}\right),
       	\een
       	where $\alpha=\alpha_1-\alpha_2-2$ and $m\in\NN_0$ such that $\alpha_3-\alpha_1+3=m$.
  	  \ele

  \proof
  Let us firstly observe that the formal power series $\hat{\psi_1}(x)$ is divergent.
  Indeed,  we have that
  $$\,
  \lim_{n\rightarrow \infty}\, F(\alpha; \beta; n)=\frac{\beta_3-\beta_2}{\beta_1-\beta_2}.
  \,$$
  Then the radius $R$ of convergence of the series $\hat{\psi_1}(x)$ is computed as
  $$\,
  R=\lim_{n\rightarrow \infty}\left|
  \frac{\beta_3-\beta_2}{2+n+\alpha_3-\alpha_2}\right|=0.
  \,$$
  
  Since when $|\beta_3-\beta_1| < |\beta_3-\beta_2|$ the number series
  $F(\alpha; \beta; n)$	are absolutely convergent for all $n\in\NN_0$  
 and the limit $ \lim_{n\rightarrow \infty}\, F(\alpha; \beta; n)$ does not depend on $n$  
 then all $|F(\alpha; \beta; n)|$ are bounded.
  
  In \leref{phi} we have proved that the series $\hat{\psi}_2(x)$ is 1-summable in every
  direction $\theta \neq \arg(\beta_2-\beta_3)$ from 0 to $+\infty \,e^{i \theta}$.
  Consider only the series $\hat{\psi}_1(x)$.
  Let $|2+\alpha_3-\alpha_2| \leq 1$. Then 
  $$\,
  |(2+\alpha_3-\alpha_2)^{(n)}| \leq n!.
  \,$$
  Let $|2+\alpha_3-\alpha_2| > 1$. Then 
  $$\,
  |(2+\alpha_3-\alpha_2)^{(n)}| <
  (|2+\alpha_3-\alpha_2|+1)^n\,n!.
  \,$$
  Therefore the formal power series
  $$\,
  \sum_{n=0}^{\infty} 
  \frac{(-1)^n\,(2+\alpha_3-\alpha_2)^{(n)}}{(\beta_3-\beta_2)^n}\,
  F(\alpha; \beta; n)\,x^n
  \,$$
  is of Gevrey order 1 with constants
  $$\,
  C=\max_{n \geq 0} |F(\alpha; \beta; n)|,\quad
  A=\frac{1}{|\beta_3-\beta_2|}
    \,$$
  if $|2+\alpha_3-\alpha_2| \leq 1$, and
  $$\,
  C=\max_{n \geq 0} |F(\alpha; \beta; n)|,  \quad
  A=\frac{|2+\alpha_3-\alpha_2|+1}{|\beta_3-\beta_2|}
  \,$$
  if $|2+\alpha_3-\alpha_2| > 1$.

  As a result the formal Borel transforms of $\hat{\psi_1}(x)$ and $\hat{\psi_2}(x)$
  \ben
  \psi_1(\xi)    &=&
  (\hat{\B}_1 \hat{\psi_1})(\xi)=
  \sum_{n=0}^{\infty}
  \frac{(-1)^n\,(2+\alpha_3-\alpha_2)^{(n)}}{(\beta_3-\beta_2)^n}\,
  \frac{F(\alpha; \beta; n)}{n!}\,\xi^n,\\[0.3ex]
  \psi_2(\xi)    &=&
  (\hat{\B}_1 \hat{\psi_2})(\xi)=
  \sum_{n=0}^{\infty}
  \frac{(-1)^n\,(2+\alpha_3-\alpha_2)^{(n)}}{(\beta_3-\beta_2)^n}\,\frac{\xi^n}{n!}\\[0.3ex]
  &=&	
  \left(1+\frac{\xi}{\beta_3-\beta_2}\right)^{\alpha_2-\alpha_3-2}
  \een
  are  analytic functions near the origin in the Borel plane. Since there the series
  $\psi_1(\xi)$ is absolutely convergent, then
  every series obtained from $\psi_1(\xi)$ by changing the positions of its terms, is also convergent
  and has the same sum. In particular, the series
  $$\,
  \psi_1(\xi)=
  \sum_{p=0}^{\infty} u_p (\xi),
  \,$$
  where 
  \ben
  u_p(\xi)  &=&
  \frac{(2+\alpha_3-\alpha_2)^{(p)}}{(4+\alpha_3-\alpha_1)^{(p)}}
  \left(\frac{\beta_3-\beta_1}{\beta_3-\beta_2}\right)^p\,
  F_p(\xi),\quad 
  p \geq 1,\\[0.2ex]
  u_0(\xi)  &=&
  \left(1+\frac{\xi}{\beta_3-\beta_2}\right)^{\alpha_2-\alpha_3-2}	
  \een
  with
  $$\,
  F_p(\xi)=\sum_{s=0}^{\infty}
  \frac{(2+p+\alpha_3-\alpha_2)^{(s)} (4+\alpha_3-\alpha_1)^{(s)}}
  {(4+p+\alpha_3-\alpha_1)^{(s)} s!}
  \left(-\frac{\xi}{\beta_3-\beta_2}\right)^s
  \,$$ 
 fulfills this property. Then  the function
  $\psi(\xi)=\psi_1(\xi) - \psi_2(\xi)$ is reduced to
  $$\,
  \psi(\xi)=\sum_{p=1}^{\infty} u_p(\xi).
  \,$$

  Ir order to see the dependence of $\psi(\xi)$ on the 
  directions $\arg (\beta_2-\beta_3)$ and $\arg (\beta_1-\beta_3)$ we use the following
  connection formula for the hypergeometric series $F(a,b; c; z)$ (see 2,10(1) in \cite{HB-AE})
   \ben
   F(a, b; c; z)  &=&
   A_1\,F(a, b; a+b-c+1; 1-z)\\[0.1ex]
                  &+&
   A_2\,(1-z)^{c-a-b}\,F(c-a, c-b; c-a-b+1; 1-z),
   \een
   where
   $$\,
    A_1=\frac{\Gamma(c)\,\Gamma(c-a-b)}{\Gamma(c-a)\,\Gamma(c-b)},\quad
    A_2=\frac{\Gamma(c)\,\Gamma(a+b-c)}{\Gamma(a)\,\Gamma(b)},
   \,$$
   provided that $|\arg(1-z)| < \pi$ and $c-a-b\notin\ZZ$. Applying this formula to the  hypergeometrc
  series $F_p(\xi)$ we transform it into
  \be\label{F1}
  F_p(\xi)   &=&
  A_1(p)\,F_p(a_1, b_1 ; c_1 ; 1+\frac{\xi}{\beta_3-\beta_2})  \\
  &+&
  A_2(p)\,\left(1+\frac{\xi}{\beta_3-\beta_2}\right)^{\alpha_2-\alpha_3-2}\,
  F_p(a_2, b_2 ; c_2 ; 1+\frac{\xi}{\beta_3-\beta_2}).\nonumber
  \ee
  The parameters are given by
  \ben
  a_1=2+p+\alpha_3-\alpha_2,
  &\quad&   a_2=2+\alpha_2-\alpha_1,\\
  b_1=4+\alpha_3-\alpha_1,
  &\quad&  b_2=p,\\
  c_1=3+\alpha_3-\alpha_2,
  &\quad& c_2=\alpha_2-\alpha_3-1,
  \een
  and the coefficients $A_1(p)$ and $A_2(p)$ are determined as
  \be\label{A}
  A_1(p)  &=&
  \frac{\Gamma(4+p+\alpha_3-\alpha_1)\,\Gamma(\alpha_2-\alpha_3-2)}
  {\Gamma(2+\alpha_2-\alpha_1)\,\Gamma(p)},\\[0.3ex]
  A_2(p)  &=&
  \frac{\Gamma(4+p+\alpha_3-\alpha_1)\,\Gamma(\alpha_3-\alpha_2+2)}
  {\Gamma(2+p+\alpha_3-\alpha_2)\,\Gamma(4+\alpha_3-\alpha_1)}\nonumber			
  \ee
  for $\left|\arg\,\left(1+\frac{\xi}{\beta_3-\beta_2}\right)\right| < \pi$.
  Under assumptions of this lemma  the parameters $a_1, b_1, b_2, c_1$ do not belong
  to $\ZZ_{\leq 0}$. The parameters $a_2$ and $c_2$ can be non-positive integers, but
  the definition of the hypergeometric series can be extended for $c\in\ZZ_{\leq 0}$
  (see the end of the proof). The coefficients $A_1(p), A_2(p)$
  are well defined provided that $\alpha_3-\alpha_2\notin\ZZ$. At the end of the proof
  we show that the coefficient $A_1(p)$ is  well determined also  when
  $\alpha_3-\alpha_2\in\ZZ_{\geq -1}$ but $\alpha_2-\alpha_1\in\ZZ_{\leq -2}$. For now
  we assume that $\alpha_3-\alpha_2\notin\ZZ$.

  Then the sum $\psi(\xi)$ is changed into
  $$\,
  \psi(\xi)=
  \sum_{p=1}^{\infty} v_{p}(\xi) +
  \sum_{p=1}^{\infty} w_{p}(\xi),
  \,$$
  where 
  \ben
  v_{p}(\xi)  &=&
  \frac{A_1(p)\,(2+\alpha_3-\alpha_2)^{(p)}}
  {(4+\alpha_3-\alpha_1)^{(p)}}
  \left(\frac{\beta_3-\beta_1}{\beta_3-\beta_2}\right)^p
  F_{p, 1}(\xi),\\[0.4ex]
  w_{p}(\xi)  &=&
  \frac{A_2(p)\,(2+\alpha_3-\alpha_2)^{(p)}}
  {(4+\alpha_3-\alpha_1)^{(p)}}
  \left(\frac{\beta_3-\beta_1}{\beta_3-\beta_2}\right)^p
  \left(1+\frac{\xi}{\beta_3-\beta_2}\right)^{\alpha_2-\alpha_3-2}\,F_{p,2}(\xi).
  	\een
  Here we denote by $F_{p, i}(\xi), \,i=1, 2$ the function
  $F_p(a_i, b_i; c_i; 1+\frac{\xi}{\beta_3-\beta_2})$.

  For any $p \geq 1$ the function $v_p(\xi)$ is an analytic function in
  $1 +\xi/(\beta_3-\beta_2)$.
  It turns out that it is not analytic in $1+\xi/(\beta_3-\beta_1)$. Indeed,
  using that $\Gamma(4+p+\alpha_3-\alpha_1)=(4+\alpha_3-\alpha_1)^{(p)}\Gamma(4+\alpha_3-\alpha_1)$,
  we rewrite $\Upsilon=\sum v_p(\xi)$ as
  $$\,
  \Upsilon=\frac{\Gamma(4+\alpha_3-\alpha_1) \Gamma(\alpha_2-\alpha_3-2)}
  {\Gamma(2+\alpha_2-\alpha_1)}
  \sum_{p=1}^{\infty}
  \frac{(2+\alpha_3-\alpha_2)^{(p)}}{\Gamma(p)}
  \left(\frac{\beta_3-\beta_1}{\beta_3-\beta_2}\right)^p
  F_{p, 1}(\xi).
  \,$$
  Next we represent this sum as a power series in $1+\xi/(\beta_3-\beta_2)$
  $$\,
  \Upsilon=\frac{\Gamma(4+\alpha_3-\alpha_1) \Gamma(\alpha_2-\alpha_3-2)}
  {\Gamma(2+\alpha_2-\alpha_1)}
  \sum_{s=0}^{\infty} C_s\,\left(1+\frac{\xi}{\beta_3-\beta_2}\right)^s.
  \,$$
  For the coefficients $C_s$ we obtain consecutively
  \ben
  C_s   &=&
  \frac{(4+\alpha_3-\alpha_1)^{(s)}}
  {(3+\alpha_3-\alpha_2)^{(s)} s!}
  \sum_{p=1}^{\infty}
  \frac{(2+\alpha_3-\alpha_2)^{(p)}\,(2+p+\alpha_3-\alpha_2)^{(s)}}{(p-1)!}
  \left(\frac{\beta_3-\beta_1}{\beta_3-\beta_2}\right)^p\\[0.25ex]
  &=&
  \frac{(4+\alpha_3-\alpha_1)^{(s)} (2+\alpha_3-\alpha_2)}{s!}
  \frac{\beta_3-\beta_1}{\beta_3-\beta_2}
  \left(\frac{\beta_1-\beta_2}{\beta_3-\beta_2}\right)
  ^{\alpha_2-\alpha_3-3-s}.
  \een
  Finally  $\Upsilon$ becomes
  \ben
  \Upsilon  &=& \sum_{p=1}^{\infty} v_p(\xi)=
  \frac{\Gamma(4+\alpha_3-\alpha_1) \Gamma(\alpha_2-\alpha_3-2) (2+\alpha_3-\alpha_2)}
  {\Gamma(2+\alpha_2-\alpha_1)}\frac{\beta_3-\beta_1}{\beta_3-\beta_2}\\[0.25ex]
  &\times&
  \left(\frac{\beta_1-\beta_2}{\beta_3-\beta_2}\right)^{\alpha_2-\alpha_3-3}
  \left(\frac{\beta_1-\beta_3}{\beta_1-\beta_2}\right)^{\alpha_1-\alpha_3-4}
  \left(1+\frac{\xi}{\beta_3-\beta_1}\right)^{\alpha_1-\alpha_3-4}.
  \een
  
  To make the computation
  of the Stokes multiplier, corresponding to the singular direction
  $\theta=\arg (\beta_2-\beta_3)$, more simple, we will represent the sum
  $\sum_{p=1}^{\infty} w_p(\xi)$  in  powers of $(1+\xi/(\beta_3-\beta_2)$
  \ben
  \sum_{p=1}^{\infty} w_p(\xi)=
  \left(1+\frac{\xi}{\beta_3-\beta_2}\right)^{\alpha_2-\alpha_3-2}
  \sum_{s=0}^{\infty} C_s\,\left(1+\frac{\xi}{\beta_3-\beta_2}\right)^s.
  \een
  For the coefficients $C_s$ we obtain consecutively
  \ben
  C_s  &=&
  \frac{(2+\alpha_2-\alpha_1)^{(s)}}{(\alpha_2-\alpha_3-1)^{(s)}\,s!}
  \sum_{p=1}^{\infty}
  (p)^{(s)}\,\left(\frac{\beta_3-\beta_1}{\beta_3-\beta_2}\right)^p\\[0.5ex]
  &=&
  \frac{(2+\alpha_2-\alpha_1)^{(s)}}{(\alpha_2-\alpha_3-1)^{(s)}}
  \frac{(\beta_3-\beta_1)\,(\beta_3-\beta_2)^s}{(\beta_1-\beta_2)^{s+1}}.
  \een
  Finally, 
  \ben
     & &
  \sum_{p=1}^{\infty}w_p(\xi)\\[0.1ex]
      &=&
  \left(1+\frac{\xi}{\beta_3-\beta_2}\right)^{\alpha_2-\alpha_3-2}
  \frac{\beta_3-\beta_1}{\beta_1-\beta_2}
  \sum_{s=0}^{\infty} 
  \frac{(2+\alpha_2-\alpha_1)^{(s)}}
  {(\alpha_2-\alpha_3-1)^{(s)}}
  \left(\frac{\beta_3-\beta_2}{\beta_1-\beta_2}\right)^s
  \left(1+\frac{\xi}{\beta_3-\beta_2}\right)^s.
  \een

  Then when $\alpha_3-\alpha_2\notin\ZZ_{\geq -1}$ for any directions $\theta_3 \neq \arg(\beta_2-\beta_3)$ from $0$ to
  $+\infty\,e^{i \theta_3}$ and $\theta_2 \neq \arg(\beta_1-\beta_3)$ from $0$ to
  $+\infty\,e^{i \theta_2}$ the associate Laplace transforms
  \ben
  \psi_{\theta}(x) &=&
  \frac{\beta_3-\beta_1}{\beta_1-\beta_2}
  \sum_{s=0}^{\infty} \frac{(2+\alpha_2-\alpha_1)^{(s)}}
  {(\alpha_2-\alpha_3-1)^{(s)}}
  \left(\frac{\beta_3-\beta_2}{\beta_1-\beta_2}\right)^s\\[0.4ex]
  &\times&
  \int_0^{+\infty e^{i \theta_3}}
  \left(1+\frac{\xi}{\beta_3-\beta_2}\right)^{\alpha_2-\alpha_3-2+s}
  e^{-\frac{\xi}{x}} d \left(\frac{\xi}{x}\right)\\[0.4ex]
  &-&
  \frac{(\beta_3-\beta_1)\Gamma(4+\alpha_3-\alpha_1)\Gamma(\alpha_2-\alpha_3-1)}
  {(\beta_3-\beta_2)\Gamma(2+\alpha_2-\alpha_1)}
  \left(\frac{\beta_1-\beta_2}{\beta_3-\beta_2}\right)^{\alpha_2-\alpha_3-3}\\[0.4ex]
  &\times&
  \left(\frac{\beta_1-\beta_3}{\beta_1-\beta_2}\right)^{\alpha_1-\alpha_3-4}     
  \int_0^{+\infty\,e^{i \theta_2}}
  \left(1+\frac{\xi}{\beta_3-\beta_1}\right)^{\alpha_1-\alpha_3-4}\,
  e^{-\frac{\xi}{x}}\,d \left(\frac{\xi}{x}\right)			
  \een
  define the corresponding 1-sum of the series $\hat{\psi}(x)$  in such  directions.
  
   We finish the proof considering the case when $\alpha_3-\alpha_2\in\ZZ_{\geq -1}$
   but $\alpha_2-\alpha_1\in\ZZ_{\leq -2}$. Since we assume that 
   $\alpha_3-\alpha_1\notin\ZZ_{\leq -4}$ then $2+\alpha_2-\alpha_1 \geq \alpha_2-\alpha_3-1$.
   Recall that when $b, c\in\ZZ_{\leq 0}$ such that $b \geq c$ the hypergeometric series
   $F(a, b; c; z)$ is reduced to a polynomial (see 2.1(4) in \cite{HB-AE})
    $$\,
       F(a, -n; -m; z)=\sum_{s=0}^n
       \frac{(a)^{(s)}\,(-n)^{(s)}}{(-m)^{(s)}\,s!}\,z^s.
    \,$$
    Observe also that when $2+\alpha_2-\alpha_1 \geq \alpha_2-\alpha_3-1$ we have
     $$\,
       \frac{\Gamma(\alpha_2-\alpha_3-1)}{\Gamma(2+\alpha_2-\alpha_1)}=
       \frac{1}{(\alpha_2-\alpha_3-1)\,(\alpha_2-\alpha_3) \dots (\alpha_2-\alpha_3+m-2)},
     \,$$
     where we denote $2+\alpha_2-\alpha_1=\alpha_2-\alpha_3-1+m,\,m\in\NN_0$.
     Combining the above two observations we transform the function $\psi_{\theta}(x)$ into
         \ben
         \psi_{\theta}(x) &=&
         \frac{\beta_3-\beta_1}{\beta_1-\beta_2}
         \sum_{s=0}^{\alpha_1-\alpha_2-2} \frac{(2+\alpha_2-\alpha_1)^{(s)}}
         {(\alpha_2-\alpha_3-1)^{(s)}}
         \left(\frac{\beta_3-\beta_2}{\beta_1-\beta_2}\right)^s\\[0.4ex]
         &\times&
         \int_0^{+\infty e^{i \theta_3}}
         \left(1+\frac{\xi}{\beta_3-\beta_2}\right)^{\alpha_2-\alpha_3-2+s}
         e^{-\frac{\xi}{x}} d \left(\frac{\xi}{x}\right)\\[0.4ex]
         &-&
         \frac{(\beta_3-\beta_1)\Gamma(4+\alpha_3-\alpha_1)}
         {(\beta_3-\beta_2)(\alpha_2-\alpha_3-1) (\alpha_2-\alpha_3) \dots (\alpha_2-\alpha_3+m-2)}
         \left(\frac{\beta_1-\beta_2}{\beta_3-\beta_2}\right)^{\alpha_2-\alpha_3-3}\\[0.4ex]
         &\times&
         \left(\frac{\beta_1-\beta_3}{\beta_1-\beta_2}\right)^{\alpha_1-\alpha_3-4}     
         \int_0^{+\infty\,e^{i \theta_2}}
         \left(1+\frac{\xi}{\beta_3-\beta_1}\right)^{\alpha_1-\alpha_3-4}\,
         e^{-\frac{\xi}{x}}\,d \left(\frac{\xi}{x}\right).			
         \een
    
    Note that when $\alpha_2-\alpha_1\in\ZZ_{\leq -2}$ but $\alpha_3-\alpha_2\notin\ZZ$
    the function $\psi_{\theta}(x)$ is reduced to
             \ben
             \psi_{\theta}(x) &=&
             \frac{\beta_3-\beta_1}{\beta_1-\beta_2}
             \sum_{s=0}^{\alpha_1-\alpha_2-2} \frac{(2+\alpha_2-\alpha_1)^{(s)}}
             {(\alpha_2-\alpha_3-1)^{(s)}}
             \left(\frac{\beta_3-\beta_2}{\beta_1-\beta_2}\right)^s\\[0.4ex]
             &\times&
             \int_0^{+\infty e^{i \theta_3}}
             \left(1+\frac{\xi}{\beta_3-\beta_2}\right)^{\alpha_2-\alpha_3-2+s}
             e^{-\frac{\xi}{x}} d \left(\frac{\xi}{x}\right).
             \een
             Therefore in this case the direction $\theta_2=\arg(\beta_1-\beta_3)$
             is not a singular direction for the function $\psi_{\theta}(x)$.
             
  This completes the proof.
  
  \qed

  \bre{r2}
  Observe that the functions $u_{p, \theta}$ can be regarded as the 1-sums
  of the formal series
  \ben
  \hat{u_p}(x)   &=&
  \frac{(2+\alpha_3-\alpha_2)_p}{(4+\alpha_3-\alpha_1)_p}
  \left(\frac{\beta_3-\beta_1}{\beta_3-\beta_2}\right)^p\\[0.3ex]
  &\times&
  \sum_{n=0}^{\infty}
  \frac{(2+p+\alpha_3-\alpha_2)^{(n)}\,(4+\alpha_3-\alpha_1)^{(n)}}
  {(4+p+\alpha_3-\alpha_1)^{(n)}}\,
  \frac{(-1)^n\,x^n}{(\beta_3-\beta_2)^n}.		
  \een
  In this representation we look on the formal series $\hat{\psi_1}(x)$ as
  \ben
  \hat{\psi_1}(x) =
  \sum_{p=0}^{\infty} \hat{u_p}(x).
  \een
  \qed
  \ere

  \bre{transf}
    The transformation formula \eqref{F1} does not remain valid when
    $\alpha_3-\alpha_2\in\ZZ_{\geq -1}$ but $\alpha_2-\alpha_1\notin\ZZ_{\leq -2}$. 
    The reason is that the Gamma function
     $\Gamma(z)$ has simple poles at $z=0, -1, -2, -3, \dots$.
     Then for $\alpha_3-\alpha_2\in\ZZ_{\geq -1}$ the connection coefficient
     $A_1(p)$ in \eqref{A} is infinite while the series $F_p(\xi)$ is convergent.
     To overcome this problem one can use the so called logarithmic cases of the
     analytic continuation of the hypergeometric series. These cases are too complicated
     for  a computation by hand and for this reason they are out of scope of this paper.
     \qed
  \ere

  \bre{1-sum-2}
  Represent $\psi_{\theta}(x)$ as
   $$\,
    \psi_{\theta}(x)=\varpi_{\theta_3}(x) - \chi_{\theta_2}(x).
   \,$$
  When we move the direction $\theta_2$ (resp. $\theta_3$) continuosly the corresponding 1-sums
  $\chi_{\theta_2}(x)$ (resp. $\varpi_{\theta_3}(x)$) stick each other analytically
  and define a holomorphic function $\tilde{\chi}(x)$ (resp. $\tilde{\varpi}(x)$)
  on a sector of opening $3 \pi$. Denote by $\theta$ the singular direction of
  $\chi_{\theta_2}(x)$ or $\varpi_{\theta_3}(x)$.
  When $\pi/2 < \theta \leq 3 \pi/2$ this sector is defined by
  $\theta-2 \pi-\pi/2 < \arg (x) <\theta + \pi/2$. When $-\pi/2 < \theta \leq \pi/2$
  this sector is defined by $\theta-\pi/2 < \arg (x) < \theta+2 \pi +\pi/2$
  (see paragraph 5.10 in \cite{Sa}).
  In every non-singular direction $\theta_2$ (resp. $\theta_3$) the multivalued function $\tilde{\chi}(x)$
  (resp. $\tilde{\varpi}(x)$) has one value $\chi_{\theta_2}(x)$ (resp.  $\varpi_{\theta_3}(x)$).
  Near the corresponding singular direction $\theta$ the functions
  $\tilde{\chi}(x)$ and $\tilde{\varpi}(x)$ have two different values:
  $\chi^+_{\theta}(x)=\chi_{\theta+\epsilon}(x),\,
  \varpi^+_{\theta}(x)=\varpi_{\theta+\epsilon}(x)$ and
  $\chi^-_{\theta}(x)=\chi_{\theta-\epsilon}(x),\,
  \varpi^-_{\theta}(x)=\varpi_{\theta-\epsilon}(x)$, where $\epsilon > 0$ is a
  small number.
  \qed
  \ere

  We have a similar  summable result for the formal series $\hat{\phi}(x)$.
  
  \ble{psi2}
  Assume that $|\beta_3-\beta_1| < |\beta_3-\beta_2|$ provided that $\beta_j$'s
  are distinct . 
  Assume also that
   $\alpha_3-\alpha_1\in\ZZ_{\leq -4}$ with $\alpha_2-\alpha_1\notin\ZZ_{\geq -1}$.
  Then the formal power series $\hat{\phi}(x)$  is 1-summable and 
  	for any directions $\theta_2 \neq \arg (\beta_1-\beta_3)$ from $0$
  	to $+\infty\,e^{i \theta_2}$ and $\theta \neq \arg(\beta_i-\beta_3),\,i=1, 2$ from $0$ to
  	$+\infty \,e^{i \theta}$    the function
  	\ben
  	\phi_{\theta}(x) &=&
  	\left(\frac{\beta_1-\beta_2}{\beta_3-\beta_2}\right)^{\alpha_2-\alpha_1+1}
  	\int_0^{+\infty e^{i \theta_2}}
  	\left(1+\frac{\xi}{\beta_3-\beta_1}\right)^{-1}\,
  	e^{-\frac{\xi}{x}}\,d \left(\frac{\xi}{x}\right)\\[0.4ex]
  	&-&
  	\int_0^{+\infty e^{i \theta}}
  	\left(1+\frac{\xi}{\beta_3-\beta_2}\right)^{1+\alpha_2-\alpha_1}
  	\left(1+\frac{\xi}{\beta_3-\beta_1}\right)^{-1}\,
  	e^{-\frac{\xi}{x}}\,d \left(\frac{\xi}{x}\right)
  	\een
  	defines the corresponding 1-sum in such  directions.
  	  \ele

  Combining \leref{phi}, \leref{psi1} and \leref{psi2} we obtain an
  actual fundamental matrix at the origin.

  \bpr{as-2}
  Assume that $\beta_j$'s are distinct such that $|\beta_3-\beta_1| < |\beta_3-\beta_2|$. 
   Then for every non-singular direction $\theta$ the  equation \eqref{initial} 
  possesses an unique actual fundamental matrix $\Phi_{\theta}(x)$ at the origin in the form
   \be\label{xyz}
     \Phi_{\theta}(x)=H_{\theta}(x)\,(x^{\Lambda}\,exp\left(-\frac{Q}{x}\right))_{\theta},
   \ee
   where the matrices $\Lambda$ and $Q$ are defined by \eqref{l-q} and 
   $(x^{\Lambda}\exp(-Q/x))_{\theta}$ is the branch of $x^{\Lambda}\exp(-Q/x)$ for $\arg(x)=\theta$.
   In particular, $\Phi_{\theta + 2 \pi}(x)=\Phi_{\theta}(x)\,\hat{M}$. 
  The matrix $H_{\theta}(x)$ is given by
  \be\label{h-2}
  H_{\theta}(x)=\left(\begin{array}{ccc}
  	1   &\frac{x^2\,\varphi_{\theta_1}(x)}{\beta_2-\beta_1}       
  	&\frac{x^4\,\psi(x)}{(\beta_3-\beta_2)(\beta_3-\beta_1)}\\[0.3ex]
  	0   &1                                                   &\frac{x^2\,\varphi_{\theta_2}(x)}{\beta_3-\beta_2}\\[0.2ex]
  	0   &0                                                   &1
  \end{array}
  \right).			
  \ee
  The actual function $\varphi_{\theta_i}(x)$ is  defined by \eqref{s2} if
  $\alpha_j-\alpha_i\in\ZZ_{\leq -2}$, and by \eqref{al} otherwise. 
  
  The function  $\psi(x)$ is defined 
  as follows:
  \begin{enumerate}
  	\item\,
  	If $\alpha_3-\alpha_1\notin\ZZ_{\leq -4}$ and $\alpha_3-\alpha_2\notin\ZZ_{\leq -2}$ then
          \ben
          \psi(x)=
          F(\alpha; \beta; 0)\,
          \int_0^{+\infty e^{i \theta_2}}
          \left(1+\frac{\nu}{\beta_3-\beta_1}\right)^{\alpha_1-\alpha_3-4}
          e^{-\frac{\nu}{x}} d \left(\frac{\nu}{x}\right)
          - \psi_{\theta}(x),  
          \een
          where $\psi_{\theta}(x)$ is defined by the first part of \leref{psi1}. The number series $F(\alpha; \beta; 0)$ 
          is defined by \eqref{F0} for $\sigma=\infty$.
  
      \item\,
         If $\alpha_3-\alpha_1\in\ZZ_{\geq -3}$ but $\alpha_2-\alpha_1\in\ZZ_{\leq -2}$ then
         $\psi(x)$ is given as above, but the function $\psi_{\theta}(x)$ is defined by the second 
         part of \leref{psi1}.
  
      \item\,
      If $\alpha_3-\alpha_1\notin\ZZ_{\leq -4}$ but $\alpha_3-\alpha_2\in\ZZ_{\leq -2}$ then
           \ben
           \psi(x)=
           F(\alpha; \beta; 0)\,
           \int_0^{+\infty e^{i \theta_2}}
           \left(1+\frac{\nu}{\beta_3-\beta_1}\right)^{\alpha_1-\alpha_3-4}
           e^{-\frac{\nu}{x}} d \left(\frac{\nu}{x}\right)
           - \sum_{n=0}^{\alpha_2-\alpha_3-2} a_n\,x^n,  
           \een
           where  $F(\alpha; \beta; 0)$ 
           is defined by \eqref{F0} for $\sigma=\alpha_2-\alpha_3-2$. The coefficients $a_n$
           are given by \eqref{an}.
  
      \item\,
         If $\alpha_3-\alpha_1\in\ZZ_{\leq -4}$ with $\alpha_3-\alpha_2\notin\ZZ$ then
            		\ben
            	\psi(x) &=& 
            		\sum_{l=0}^{\alpha_1-\alpha_3-4}
            		b_l\,x^l\left(\sum_{s=0}^{\alpha_1-\alpha_3-4-l}
            		(-1)^s \frac{(4+l+\alpha_3-\alpha_1)^{(s)}}{(\beta_3-\beta_1)^s}\,x^s\right)\\[0.4ex]
            		&+&
            		x^{\alpha_1-\alpha_3-3} \,b_{\alpha_1-\alpha_3-3}
                  		\int_0^{+\infty e^{i \theta}}
            		\left(1+\frac{\xi}{\beta_3-\beta_2}\right)^{1+\alpha_2-\alpha_1}
            		\left(1+\frac{\xi}{\beta_3-\beta_1}\right)^{-1}
            		e^{-\frac{\xi}{x}} d \left(\frac{\xi}{x}\right),     
            	\een	  
            	where
           $$\,
             b_n=(-1)^n\,\frac{(2+\alpha_3-\alpha_2)^{(n)}}{(\beta_3-\beta_2)^n}.
           \,$$ 	
  
       \item\,
            If $\alpha_3-\alpha_1\in\ZZ_{\leq -4}$ and $\alpha_2-\alpha_1\in\ZZ_{\geq -1}$ then
                 	\ben
                 	\psi(x) &=&
                 	\sum_{l=0}^{\alpha_1-\alpha_3-4}
                 	b_l\,x^l\left(\sum_{s=0}^{\alpha_1-\alpha_3-4-l}
                 	(-1)^s \frac{(4+l+\alpha_3-\alpha_1)^{(s)}}{(\beta_3-\beta_1)^s}\,x^s\right)\\[0.4ex]
                 	&+&
                 	x^{\alpha_1-\alpha_3-3} \,b_{\alpha_1-\alpha_3-3}
                 	\int_0^{+\infty e^{i \theta}}
                 	\left(1+\frac{\nu}{\beta_3-\beta_1}\right)^{-1}
                 	e^{-\frac{\nu}{x}} d \left(\frac{\nu}{x}\right)\\[0.4ex]
                 	&-&
                 	x^{\alpha_1-\alpha_3-3} \,b_{\alpha_1-\alpha_3-3}
                 	\sum_{n=0}^{\alpha_2-\alpha_1} a_n\,x^n,  
                 	\een
                 where $b_n$ are defined as above. The coefficients $a_n$ are given
                 by \eqref{ant}. The series $\sum_{n=0}^{\alpha_2-\alpha_1} a_n\,x^n$
                 is equal to zero if $\alpha_2-\alpha_1=-1$.
  
            \item\,
            If $\alpha_2-\alpha_1,\,\alpha_3-\alpha_2\in\ZZ_{\leq -2}$ then                 
                       $$\,
                       \psi(x)=
                       \sum_{l=0}^{\alpha_2-\alpha_3-2} b_l\,x^l
                       \left(\sum_{s=0}^{\alpha_1-\alpha_3-4-l} (-1)^s\,
                       \frac{(4+l+\alpha_3-\alpha_1)^{(s)}}{(\beta_3-\beta_1)^s}\,x^s\right).
                       \,$$              
  \end{enumerate}
      For a singular direction $\theta$ the equation \eqref{initial} has two actual fundamental
      matrix at the origin
       $$\,
        \Phi^{\pm}_{\theta}(x)=\Phi_{\theta \pm \epsilon}(x),
       \,$$
       where the matrices $\Phi_{\theta \pm \epsilon}$ are given by \eqref{xyz} for a small 
       positive number $\epsilon$.
       
      Moreover the matrix $\Phi_{\theta}(x)$ given by the above formulas and the fundamental
      matrix $\Phi(x)$ introduced by equations \eqref{global-m}-\eqref{global-el}
      define the same actual fundamental matrix at the origin.
 \epr

    \proof

    From \thref{actual}, \prref{fs-0},  \leref{phi}, \leref{psi1}, \leref{psi2}
    and \reref{1-sum-1} it follows that the matrix $\Phi_{\theta}(x)$ defined as above is
    the unique actual fundamental matrix at the origin, associate with the formal
    fundamental matrix $\hat{\Phi}(x)$.
    
    Since the element $\Phi_{13}(x)$ is obtained by the iterated integral
     $$\,
            \Phi_{13}(x) =
            \Phi_1(x)\int_{\gamma_3(x)}
            \frac{\Phi_2(t)}{\Phi_1(t)}
            \left( \int_{\gamma_2(t)}
            \frac{\Phi_3(t_1)}{\Phi_2(t_1)} \d t_1\right)\d t          
     \,$$
     we have only to show that the elements $\Phi_{12}(x)$ and $\Phi_{23}(x)$
     introduced by \eqref{global-el} coincide with these ones given by the 
     present proposition.

    The element $\Phi_{12}(x)$ from \eqref{global-el} can be expressed as
     \ben
      \Phi_{12}(x)
          &=&
          \Phi_1(x)\int_0^x \frac{\Phi_2(t)}{\Phi_1(t)}\, d t=
      x^{\alpha_1} e^{-\frac{\beta_1}{x}}
      \int_0^x t^{\alpha_2-\alpha_1} e^{-\frac{\beta_2-\beta_1}{t}}\,d t\\[0.4ex]
         &=&
        x^{\alpha_1} e^{-\frac{\beta_2}{x}}
        \int_0^x t^{\alpha_2-\alpha_1} e^{-\frac{\beta_2-\beta_1}{t}}
        e^{\frac{\beta_2-\beta_1}{x}}\, d t, 
     \een
     where the path $\gamma_1(x)$ is a path from 0 to $x$ approaching 0
     in the direction $\theta=\arg (\beta_2-\beta_1)$.
     Introducing a new variable $\zeta$ via
      $$\,
        (\beta_2-\beta_1)\left(\frac{1}{x} - \frac{1}{t}\right)=
        -\frac{\zeta}{x}
      \,$$
      we get
      $$\,
        \Phi_{12}(x)=\frac{x^{\alpha_2+2} e^{-\frac{\beta_2}{x}}}{\beta_2-\beta_1}
        \int_0^{+\infty}\left(
        1+\frac{\zeta}{\beta_2-\beta_1}\right)^{\alpha_1-\alpha_2-2}
        e^{-\frac{\zeta}{x}}\, d \left(\frac{\zeta}{x}\right).
      \,$$

     Similarly, the function $\Phi_{23}(x)$ from \eqref{global-el} can be
     expressed as
       $$\,
          \Phi_{23}(x)=x^{\alpha_2} e^{-\frac{\beta_3}{x}}
          \int_0^x t^{\alpha_3-\alpha_2}\,e^{-\frac{\beta_3-\beta_2}{t}}\,
          e^{\frac{\beta_3-\beta_2}{x}}\, d t
       \,$$
          where the path $\gamma_2(x)$ is a path from 0 to $x$ approaching 0
          in the direction $\theta=\arg (\beta_3-\beta_2)$.
          After  introducing a new variable $\varsigma$ via
           $$\,
           (\beta_3-\beta_2)\left(\frac{1}{x} - \frac{1}{t}\right)=
           -\frac{\varsigma}{x}
           \,$$
           the function $\Phi_{23}(x)$ is transformed to
           $$\,
           \Phi_{23}(x)=\frac{x^{\alpha_3+2} e^{-\frac{\beta_3}{x}}}{\beta_3-\beta_2}
           \int_0^{+\infty}\left(
           1+\frac{\varsigma}{\beta_3-\beta_2}\right)^{\alpha_2-\alpha_3-2}
           e^{-\frac{\varsigma}{x}}\, d \left(\frac{\varsigma}{x}\right).
           \,$$
    Analytic continuations of the so constructed $\Phi_{12}(x)$ and $\Phi_{23}(x)$ on 
    $x$-plane yield analytic functions 
    $$\,
      (\Phi_{12}(x))_{\theta_1}=\frac{x^{\alpha_2+2}\,e^{-\frac{\beta_2}{x}}}{\beta_2-\beta_1}\,
      \varphi_{\theta_1}(x),\quad
         (\Phi_{23}(x))_{\theta_2}=\frac{x^{\alpha_3+2}\,e^{-\frac{\beta_3}{x}}}{\beta_3-\beta_2}\,
         \varphi_{\theta_2}(x)
    \,$$ 
    on all non-singular directions $\theta_1$ and $\theta_2$. Here $\varphi_{\theta_j},\,j=1, 2$
    are defined by \eqref{al}.
    
    This ends the proof. 
    \qed

   It turns out that the above results remain valid for all values of the parameters $\beta_j$'s
   for which $|\beta_3-\beta_1|=|\beta_3-\beta_2|$ but $\beta_3-\beta_1\neq \pm (\beta_3-\beta_2)$
   and $\R (\alpha_2-\alpha_1) > -1$. Recall that when 
    $|\beta_3-\beta_1|=|\beta_3-\beta_2|$ but 
     $\R (\alpha_2-\alpha_1) > -1$ the number series $F$ and $\tilde{F}$ are still absolutely convergent. 
    We avoid the very particular case when $\beta_3-\beta_1=\beta_3-\beta_2$ since it implies
    that $\beta_2=\beta_1$. The case $\beta_3-\beta_1=-(\beta_3-\beta_2)$ is out of scope of this paper
    since $F(a, b; c; -1)$ has different values depending on the parameters $a, b$ and $c$.
    Then

     \bpr{as-3}
       Assume that $\beta_j$'s are distinct such that $|\beta_3-\beta_1|=|\beta_3-\beta_2|$
       but $\beta_3-\beta_1 \neq \pm (\beta_3-\beta_2)$. Assume also that $\R (\alpha_2-\alpha_1) > -1$. 
       Then an actual fundamental matrix at the origin of  equation \eqref{initial} is defined by
       \prref{as-2} with the exception of the items $(2)$ and $(6)$ provided that
       $\Re (\alpha_2-\alpha_1) > -1$. The item $(5)$ remains valid for $\alpha_3-\alpha_1\in\ZZ_{\leq -4}$
       but $\alpha_2-\alpha_1 \in\ZZ_{\geq 0}$.
    \epr


		\section{Stokes matrices}

   In this section we compute the Stokes matrices with respect to the 
   formal and actual fundamental matrices constructed in the previous two sections.

    Next example is a model of the jump  of the functions including
		 in the actual matrix  solution of the  equation \eqref{initial}.
		 
		\bex{e1}
      Consider the function
			  \be\label{f}
				   f_{\theta}(x)=\int_0^{+\infty\,e^{i \theta}}
					\left(1 + \frac{\eta}{\beta_j-\beta_i}\right)^{\alpha_i-\alpha_j-2}\,
					e^{-\frac{\eta}{x}}\, d \left(\frac{\eta}{x}\right),
				\ee
				where the path of integration is taken in the direction $\theta \neq d=\arg( \beta_i-\beta_j)$
				from $0$ to $+\infty\,e^{i \theta}$. Here $\alpha_i-\alpha_j, \beta_i-\beta_j\in\CC$ such that
			 $\beta_j$'s are distinct and $\alpha_j-\alpha_i \notin\ZZ_{\leq -2}$.
				
				The function $f_{\theta}(x)$ is an analytic function in all non-singular directions
				$\theta \neq d=\arg( \beta_i-\beta_j)$. Let $\theta^+=d + \epsilon$ and
				$\theta^-=d - \epsilon$, where $\epsilon > 0$ is a small number,
				be two non-singular neighboring directions of the singular direction $d=\arg ( \beta_i-\beta_j)$.
				Denote by $f^{+}_{\theta}(x)$ and $f^{-}_{\theta}(x)$  the integrals \eqref{f} in such
				direction respectively. Comparing these two functions, we have that
				  $$\,
					  f^{-}_{\theta}(x)=f^{+}_{\theta}(x) +
						  x^{-1}\,\int_\gamma
							 \left(1 + \frac{\eta}{\beta_j-\beta_i}\right)^{\alpha_i-\alpha_j-2}\,e^{-\frac{\eta}{x}}\,
							  d \eta,
					\,$$
					where $\gamma=\theta^{-} - \theta^+$. Without changing the integral, we can deform
					$\gamma$ in a Hankel type path $\gamma'$, going along the direction $d$
					from infinity to $\beta_i-\beta_j$, encircle $ \beta_i-\beta_j$ in the positive direction and back to
					infinity in the positive sense. Then, since 
					$\arg\left(1+\frac{\eta}{\beta_j-\beta_i}\right)=-\pi$ when 
					$\pi/2 < \arg (\beta_j-\beta_i) \leq 3 \pi/2$, the integral on $\gamma'$ becomes
					 \ben
					     & &
					  x^{-1}(-e^{- \pi\,i\,(\alpha_i-\alpha_j-2)} + e^{\pi\,i\,(\alpha_i-\alpha_j-2)})
						\int_{\beta_i-\beta_j}^{+\infty\,e^{i d}}
						\left(\frac{\eta}{\beta_i-\beta_j}-1\right)^{\alpha_i-\alpha_j-2}\,e^{-\frac{\eta}{x}}\,
						 d \eta\\[0.4ex]
						   &=&
						(\beta_i-\beta_j)\,
						 x^{-1}(e^{ \pi\,i\,(\alpha_i-\alpha_j-2)} - e^{-\pi\,i\,(\alpha_i-\alpha_j-2)})\,
						e^{-\frac{\beta_i-\beta_j}{x}}
						\int_0^{+\infty} u^{\alpha_i-\alpha_j-2} e^{-\frac{(\beta_i-\beta_j)\,u}{x}}\,d u\\[0.4ex]
						   & =&
							\frac{e^{ \pi\,i\,(\alpha_i-\alpha_j-2)} - e^{-\pi\,i\,(\alpha_i-\alpha_j-2)}}
							{(\beta_i-\beta_j)^{\alpha_i-\alpha_j-2}}\,
							x^{\alpha_i-\alpha_j-2}\,e^{-\frac{\beta_i-\beta_j}{x}}
							\int_0^{+\infty} \tau^{\alpha_i-\alpha_j-2}\,e^{-\tau}\, d \tau\\[0.3ex]
							&=&
							\frac{e^{ \pi\,i\,(\alpha_i-\alpha_j-2)} - e^{-\pi\,i\,(\alpha_i-\alpha_j-2)}}
							{(\beta_i-\beta_j)^{\alpha_i-\alpha_j-2}}\,
							x^{\alpha_i-\alpha_j-2}\,e^{-\frac{\beta_i-\beta_j}{x}}\,
							\Gamma(\alpha_i-\alpha_j-1)\\[0.3ex]
							&=&
							-\frac{2 \pi\,i}{(\beta_i-\beta_j)^{\alpha_i-\alpha_j-2}\,\Gamma(2-\alpha_i+\alpha_j)}\,
							x^{\alpha_i-\alpha_j-2}\,e^{-\frac{\beta_i-\beta_j}{x}}.
					\een
					When $-\pi/2 < \arg (\beta_j-\beta_i) \leq \pi/2$ the integral on $\gamma'$ becomes
					   \ben
						  & &
					x^{-1}(1-e^{2 \pi\,i(\alpha_i-\alpha_j-2)})
					\int_{+\infty\,e^{i\, d}}^{\beta_i-\beta_j}
					\left(1+\frac{\eta}{\beta_j-\beta_i}\right)^{\alpha_i-\alpha_j-2}\,
					e^{-\frac{\eta}{x}}\,d \eta\\[0.4ex]
					   &=&
					(\beta_j-\beta_i)x^{-1}(1-e^{2 \pi\,i(\alpha_i-\alpha_j-2)})
					e^{\frac{\beta_j-\beta_i}{x}}
					\int_{-\infty}^0 u^{\alpha_i-\alpha_j-2}\,e^{-\frac{(\beta_j-\beta_i)\,u}{x}}\, d u\\[0.4ex]
					  &=&
				\frac{e^{2 \pi\,i(\alpha_i-\alpha_j-2)}-1}{(\beta_j-\beta_i)^{\alpha_i-\alpha_j-2}}\,
				  x^{\alpha_i-\alpha_j-2}\,e^{\frac{\beta_j-\beta_i}{x}}\,
					\int_0^{+\infty} \tau^{\alpha_i-\alpha_j-2}\,e^{-\tau}\,d \tau\\[0.4ex]
					  &=&
							\frac{e^{2 \pi\,i(\alpha_i-\alpha_j-2)}-1}{(\beta_j-\beta_i)^{\alpha_i-\alpha_j-2}}\,
				  x^{\alpha_i-\alpha_j-2}\,e^{\frac{\beta_j-\beta_i}{x}}\,
					\Gamma(\alpha_i-\alpha_j-1)\\[0.4ex]
					 &=&
						-\frac{2 \pi\,i\,e^{\pi\,i(\alpha_i-\alpha_j)}}
						{(\beta_j-\beta_i)^{\alpha_i-\alpha_j-2}\,\Gamma(2-\alpha_i+\alpha_j)}\,
							x^{\alpha_i-\alpha_j-2}\,e^{-\frac{\beta_i-\beta_j}{x}}.
						\een
		We have used that
		 $$\,
		    \Gamma(z)\,\Gamma(1-z)=\frac{\pi}{\sin(\pi\,z)}.
		\,$$
		\eex
     At the first sight the above two jumps of the function $f_{\theta}(x)$ give
     the same results whatever the $\arg (\beta_j-\beta_i)$ is. But it is not true.
     Indeed, let $\beta_j-\beta_i=-1$ and $\alpha_i-\alpha_j-2=-1$, i.e.
       $$\,
         f_{\theta}(x)=\int_0^{+\infty e^{i \theta}}
         \frac{e^{-\frac{\eta}{x}}}
         {1-\eta}\,d \left(\frac{\eta}{x}\right).
       \,$$
       Then according to the first case the jump of $f_{\theta}(x)$ when $\theta$ crosses
       the singular direction $\RR_{+}$ is
        $$\,
          f^-_{\theta}(x)=f^+_{\theta}(x) -2 \pi\,i\,x^{-1}\,e^{-\frac{1}{x}}.
        \,$$
        Let now $\beta_j-\beta_i=1$ and $\alpha_i-\alpha_j-2=-1$, i.e.
             $$\,
             f_{\theta}(x)=\int_0^{+\infty e^{i \theta}}
             \frac{e^{-\frac{\eta}{x}}}
             {1+\eta}\,d \left(\frac{\eta}{x}\right).
             \,$$
             Then the jump is defined as
                             $$\,
               f^-_{\theta}(x)=f^+_{\theta}(x) + 2 \pi\,i\,x^{-1}\,e^{\frac{1}{x}}.
               \,$$
               Note that not only the exponential functions are different, but also the
               Stokes multipliers are different. 
		\qed

  Finally, we compute the Stokes matrices relevant to the actual fundamental solution,
  introduced by \prref{as-2}.

   \bth{Stokes-2}
   Assume that $\beta_j$'s are distinct such that $|\beta_3-\beta_1| < |\beta_3-\beta_2|$.
   Assume also that $\arg (\beta_1-\beta_2) \neq \arg (\beta_1-\beta_3) \neq \arg (\beta_2-\beta_3)
   \neq \arg (\beta_1-\beta_2)$.
      Then with
   respect to the actual fundamental matrix $\Phi_{\theta}(x)$ at the origin,
   given by \prref{as-2}, the  equation \eqref{initial} has 
   at most three singular directions $\theta_1=\arg (\beta_1-\beta_2), \theta_2=\arg (\beta_1-\beta_3)$ and
   $\theta_3=\arg (\beta_2-\beta_3)$. The corresponding Stokes matrices are given
   by
   \ben
      St_{\theta_1}=\left(\begin{array}{ccc}
      	1    &\mu_1    &0\\
      	0    &1             &0\\
      	0    &0             &1
      \end{array}
      \right),\quad
   St_{\theta_3}=\left(\begin{array}{ccc}
   	1    &0    &0\\
   	0    &1             &\mu_3\\
   	0    &0             &1
   \end{array}
   \right),\quad
   St_{\theta_2}=\left(\begin{array}{ccc}
   	1    &0             &\mu_2\\
   	0    &1             &0\\
   	0    &0             &1
   \end{array}
   \right).
   \een
   The Sokes multipliers $\mu_j,\,j=1,  3$ are defined as follows:
   \ben
    \mu_1  &=&
    \frac{2 \pi\,i}
    {(\beta_1-\beta_2)^{\alpha_1-\alpha_2-1}\,
    	\Gamma(2 - \alpha_1 + \alpha_2)}\quad 
    \textrm{if}\quad 
    \frac{\pi}{2} < \arg (\beta_2-\beta_1) \leq \frac{3 \pi}{2} ,\\[0.5ex]
        \mu_1 
        &=& 
        - \frac{2 \pi\,i\,e^{ \pi\,i(\alpha_1-\alpha_2)}}
    {(\beta_2-\beta_1)^{\alpha_1-\alpha_2-1}\,
    	\Gamma(2 - \alpha_1 + \alpha_2)}\quad 
    \textrm{if}\quad 
    -\frac{\pi}{2} < \arg (\beta_2-\beta_1) \leq \frac{ \pi}{2} ;\\[0.85ex]
   \mu_3  
      &=&
      \frac{2 \pi\,i}
   {(\beta_2-\beta_3)^{\alpha_2-\alpha_3-1}\,
   	\Gamma(2 - \alpha_2 + \alpha_3)}\quad 
   \textrm{if}\quad 
   \frac{\pi}{2} < \arg (\beta_3-\beta_2) \leq \frac{3 \pi}{2} ,\\[0.5ex]
   \mu_3 
       &=&
       - \frac{2 \pi\,i\,e^{ \pi\,i(\alpha_2-\alpha_3)}}
   {(\beta_3-\beta_2)^{\alpha_2-\alpha_3-1}\,
   	\Gamma(2 - \alpha_2 + \alpha_3)}\quad 
   \textrm{if}\quad 
   -\frac{\pi}{2} < \arg (\beta_3-\beta_2) \leq \frac{ \pi}{2}.
   \een
   The multiplier $\mu_2$ is defined as follows:
    \begin{enumerate}
    \item\,
       If $\alpha_3-\alpha_1\notin\ZZ_{\leq -4}$ then	 
   \ben
   \mu_2
       & =&
        2 \pi\,i\Big[\frac{F(\alpha; \beta; 0)}
        {(\beta_3-\beta_2)(\beta_1-\beta_3)^{\alpha_1-\alpha_3-3} \Gamma(\alpha_3-\alpha_1+4)}\\[0.2ex]
        &-&
        \frac{\Gamma(\alpha_2-\alpha_3-1)}
        {(\beta_3-\beta_2)^{\alpha_2-\alpha_3-1} (\beta_1-\beta_2)^{\alpha_1-\alpha_2-1}
        	\Gamma(2-\alpha_1+\alpha_2)}\Big] 
    \een  
   for $\alpha_3-\alpha_2\notin\ZZ$, 
    $$\,
        \mu_2=
        \frac{2 \pi\,i\,F(\alpha; \beta; 0)}
        {(\beta_3-\beta_2)(\beta_1-\beta_3)^{\alpha_1-\alpha_3-3} \Gamma(\alpha_3-\alpha_1+4)}
    \,$$
   for $\alpha_3-\alpha_2\in\ZZ_{\leq -2}$, and
      \ben
         \mu_2
         & =&
         2 \pi\,i\Big[\frac{F(\alpha; \beta; 0)}
         {(\beta_3-\beta_2)(\beta_1-\beta_3)^{\alpha_1-\alpha_3-3} \Gamma(\alpha_3-\alpha_1+4)}\\[0.2ex]
         &-&
         \frac{1}
         {(\beta_3-\beta_2)^{\alpha_2-\alpha_3-1} (\beta_1-\beta_2)^{\alpha_1-\alpha_2-1}
         (\alpha_2-\alpha_3-1) (\alpha_2-\alpha_3)\dots (\alpha_2-\alpha_3+m-2)}\Big] 
         \een  
      for $\alpha_3-\alpha_2\in\ZZ_{\geq -1}$ but $\alpha_2-\alpha_1\in\ZZ_{\leq -2}$, 
      if $\frac{\pi}{2} < \arg (\beta_3-\beta_1) \leq \frac{3 \pi}{2}$.

       If $-\frac{\pi}{2} < \arg (\beta_3-\beta_1) \leq \frac{\pi}{2}$ then
   \ben
   \mu_2 &=& -
        2 \pi\,i\,e^{\pi\,i (\alpha_1-\alpha_3)}\Big[\frac{F(\alpha; \beta; 0)}
        {(\beta_3-\beta_2)(\beta_3-\beta_1)^{\alpha_1-\alpha_3-3} \Gamma(\alpha_3-\alpha_1+4)}\\[0.2ex]
        &+&
        \frac{e^{\pi\,i (\alpha_1-\alpha_3)}\Gamma(\alpha_2-\alpha_3-1)}
        {(\beta_3-\beta_2)^{\alpha_2-\alpha_3-1} (\beta_1-\beta_2)^{\alpha_1-\alpha_2-1}
        	\Gamma(2-\alpha_1+\alpha_2)}\Big] 
         \een  
   for  $\alpha_3-\alpha_2\notin\ZZ$, 
    $$\,
       \mu_2 = -
       \frac{2 \pi\,i\,e^{\pi\,i (\alpha_1-\alpha_3)}\,F(\alpha; \beta; 0)}
       {(\beta_3-\beta_2)(\beta_3-\beta_1)^{\alpha_1-\alpha_3-3} \Gamma(\alpha_3-\alpha_1+4)}
    \,$$
   for $\alpha_3-\alpha_2\in\ZZ_{\leq -2}$,
   and
     \ben
    \mu_2 &=& -
    2 \pi\,i\,e^{\pi\,i (\alpha_1-\alpha_3)}\Big[\frac{F(\alpha; \beta; 0)}
    {(\beta_3-\beta_2)(\beta_3-\beta_1)^{\alpha_1-\alpha_3-3} \Gamma(\alpha_3-\alpha_1+4)}\\[0.2ex]
    &+&
    \frac{e^{\pi\,i (\alpha_1-\alpha_3)}}
    {(\beta_3-\beta_2)^{\alpha_2-\alpha_3-1} (\beta_1-\beta_2)^{\alpha_1-\alpha_2-1}
    (\alpha_2-\alpha_3-1)(\alpha_2-\alpha_3)\dots (\alpha_2-\alpha_3+m-2)}\Big] 
     \een
    for $\alpha_3-\alpha_2\in\ZZ_{\geq -1}$ but $\alpha_2-\alpha_1\in\ZZ_{\leq -2}$.  
       The number series $F(\alpha; \beta; 0)$ is given by \eqref{F0},
       as $\sigma=\alpha_2-\alpha_3-2$ if $\alpha_3-\alpha_2\in\ZZ_{\leq -2}$
       and $\sigma=\infty$ otherwise. Here $m\in\NN_0$ such that $\alpha_3-\alpha_1+3=m$.  
   
   \item\,
      If $\alpha_3-\alpha_1\in\ZZ_{\leq -4}$ then
     $$\,
      \mu_2=-\frac{ 2 \pi\,i\,\Gamma(\alpha_2-\alpha_3-1)}
      {(\beta_3-\beta_2)^{\alpha_1-\alpha_3-2}\,\Gamma(2+\alpha_2-\alpha_1)},
     \,$$
     for $\alpha_2-\alpha_1\in\ZZ_{\geq -1}$, and
     $$\,
      \mu_2=-\frac{2 \pi\,i\,(-1)^{\alpha_1-\alpha_3-3}\,\Gamma(\alpha_1-\alpha_2-1)}
        {(\beta_3-\beta_2)^{\alpha_2-\alpha_3-1}\,(\beta_1-\beta_2)^{\alpha_1-\alpha_2-1}\,
        	\Gamma(2+\alpha_3-\alpha_2)}
     \,$$
     for $\alpha_3-\alpha_2\notin\ZZ$
      and if
     $\frac{\pi}{2} < \arg (\beta_3-\beta_1) \leq \frac{3 \pi}{2}$.
     
        If $-\frac{\pi}{2} < \arg (\beta_3-\beta_1) \leq \frac{\pi}{2}$ then
           $$\,
           \mu_2=\frac{ 2 \pi\,i\,\Gamma(\alpha_2-\alpha_3-1)}
           {(\beta_3-\beta_2)^{\alpha_1-\alpha_3-2}\,\Gamma(2+\alpha_2-\alpha_1)},
           \,$$
           for $\alpha_2-\alpha_1\in\ZZ_{\geq -1}$, and
           $$\,
           \mu_2=\frac{2 \pi\,i\,(-1)^{\alpha_1-\alpha_3-3}\,\Gamma(\alpha_1-\alpha_2-1)}
           {(\beta_3-\beta_2)^{\alpha_2-\alpha_3-1}\,(\beta_1-\beta_2)^{\alpha_1-\alpha_2-1}\,
           	\Gamma(2+\alpha_3-\alpha_2)}
           \,$$
           for $\alpha_3-\alpha_2\notin\ZZ$.
   \end{enumerate}
  \ethe

  \proof

  The computation of the Stokes multipliers $\mu_j\,j=1, 3$ with respect to the functions
  $\Phi_{12}(x),\,\Phi_{23}(x)$, respectively, is a direct application of
  the definition and \exref{e1}. So, we omit their derivation as an easy exercise.
   
   We consider in details the computation of $\mu_2$ and the verification of $\mu_3$ with
   respect to $\Phi_{13}(x)$ when $\alpha_3-\alpha_1\in\ZZ_{\leq -4}$
   but $\alpha_3-\alpha_2\notin\ZZ$. In this case $\Phi_{13}(x)$ has the form
           		\ben
           		& &
           		\frac{x^{\alpha_3+4}\,e^{-\frac{\beta_3}{x}}}
           			{(\beta_3-\beta_2) (\beta_3-\beta_1)} 
           		\sum_{l=0}^{\alpha_1-\alpha_3-4}
           		b_l\,x^l\left(\sum_{s=0}^{\alpha_1-\alpha_3-4-l}
           		(-1)^s \frac{(4+l+\alpha_3-\alpha_1)^{(s)}}{(\beta_3-\beta_1)^s}\,x^s\right)\\[0.4ex]
           		&+&
           		\frac{x^{\alpha_1+1}\,e^{-\frac{\beta_3}{x}} \,b_{\alpha}}
           		 	{(\beta_3-\beta_2) (\beta_3-\beta_1)} 
           		\int_0^{+\infty e^{i \theta}}
           		\left(1+\frac{\xi}{\beta_3-\beta_2}\right)^{1+\alpha_2-\alpha_1}
           		\left(1+\frac{\xi}{\beta_3-\beta_1}\right)^{-1}
           		e^{-\frac{\xi}{x}} d \left(\frac{\xi}{x}\right),     
           		\een	  
           		where
           		$$\,
           		\alpha=\alpha_1-\alpha_3-3,\quad
           		b_n=(-1)^n\,\frac{(2+\alpha_3-\alpha_2)^{(n)}}{(\beta_3-\beta_2)^n}.
           		\,$$ 	
  Under assumptions of the theorem the function $\Phi_{13}(x)$ has two singular directions  
    $\theta_3=\arg (\beta_2-\beta_3)$  and $\theta_2=\arg (\beta_1-\beta_3)$ provided that.
     Let $\theta^+_3=\theta_3+\epsilon$ and $\theta^-_3=\theta_3-\epsilon$, where
     $\epsilon > 0$ is a small number, be two neighboring directions of the singular
     direction $\theta_3$. Let $\Phi^+_{13}(x)$ and
     $\Phi^-_{13}(x)$ be the associated  actual functions to
     these directions. Let us firstly express the term for wich $\theta_3$
     is a singular direction as
     $$\,
     	\frac{x^{\alpha_1+1}\,e^{-\frac{\beta_3}{x}} \,b_{\alpha_1-\alpha_3-3}}
     	{(\beta_3-\beta_2) (\beta_2-\beta_1)} 
     	\sum_{k=0}^{\infty} (-1)^k
     	\left(\frac{\beta_3-\beta_2}{\beta_2-\beta_1}\right)^k
     	\int_0^{+\infty e^{i \theta}}
     	\left(1+\frac{\xi}{\beta_3-\beta_2}\right)^{1+\alpha_2-\alpha_1+k}
     	e^{-\frac{\xi}{x}} d \left(\frac{\xi}{x}\right).
     \,$$ 
      Then comparing the functions $\Phi^+_{13}(x)$ and $\Phi^-_{13}(x)$
     and applying \exref{e1} we find that
      \ben
        \Phi^-_{13}(x) &=&
        \Phi^+_{13}(x) \\[0.1ex]
                          &+& 
        \frac{2 \pi\,i\,b_{\alpha_1-\alpha_3-3}}
        {(\beta_2-\beta_3)^{2+\alpha_2-\alpha_1} \Gamma(-1-\alpha_2+\alpha_1)}
        \frac{x^{\alpha_2+2}\,e^{-\frac{\beta_2}{x}}}{\beta_2-\beta_1}
         \sum_{s=0}^{\infty} 
         \frac{(-1)^s (2+\alpha_2-\alpha_1)^{(s)}}{(\beta_2-\beta_1)^s}\,x^s,
       \een  
       where we have used that
       $$\,
        \frac{1}{\Gamma(-1-\alpha_2+\alpha_1-k)}=
        \frac{(-1)^k\,(2+\alpha_2-\alpha_1)^{(k)}}{\Gamma(-1-\alpha_2+\alpha_1)}.
       \,$$
       Finally, from
       $$\,
        b_{\alpha_1-\alpha_3-3}=(-1)^{\alpha_1-\alpha_3-3}
        \frac{(2+\alpha_3-\alpha_2)^{(\alpha_1-\alpha_3-3)}}{(\beta_3-\beta_2)^{\alpha_1-\alpha_3-3}}=
        \frac{\Gamma(-1-\alpha_2+\alpha_1)}{\Gamma(2+\alpha_3-\alpha_2)\,(\beta_2-\beta_3)^{\alpha_1-\alpha_3-3}}
       \,$$ 
       we obtain that
       $$\,
         \Phi^-_{13}(x)=
                \Phi^+_{13}(x)
          + \mu_3\,\Phi_{12}(x),
      \,$$
  if    $\frac{\pi}{2} < \arg (\beta_3-\beta_2) \leq \frac{3 \pi}{2}$.
  If    $-\frac{\pi}{2} < \arg (\beta_3-\beta_2) \leq \frac{ \pi}{2}$ we find that
          \ben
          \Phi^-_{13}(x)=
          \Phi^+_{13}(x) -
         \frac{2 \pi\,i\,e^{ \pi\,i(\alpha_2-\alpha_3)}}
             {(\beta_3-\beta_2)^{\alpha_2-\alpha_3-1}\,
             	\Gamma(2 - \alpha_2 + \alpha_3)}\,\Phi_{12}(x).
          \een 
         Thus when $\alpha_3-\alpha_1\in\ZZ_{\leq -4}$ but $\alpha_3-\alpha_2\notin\ZZ$
         the mulpiplier $\mu_3$ coming from $\Phi_{13}(x)$ coincides with $\mu_3$
         coming from $\Phi_{23}(x)$.
       
       The direction $\theta_2=\arg (\beta_1-\beta_3)$ is a singular direction
       only to the solution  $\Phi_{13}(x)$.
       Let $\theta^+_2=\theta_2+\epsilon$ and $\theta^-_2=\theta_2-\epsilon$, where
       $\epsilon > 0$ is a small number, be two neighboring directions of the singular
       direction $\theta_2$. Let   $\Phi^+_{13}(x)$ and $\Phi^-_{13}(x)$ be the associated  actual 
       solutions to these directions. Then comparing the functions $\Phi^+_{13}(x)$ and $\Phi^-_{13}(x)$
        we find that
             \ben
             \Phi^-_{13}(x) &=&
             \Phi^+_{13}(x) \\[0.1ex]
             &-&
             \frac{2 \pi\,i (-1)^{\alpha_1-\alpha_3-3}\,\Gamma(\alpha_1-\alpha_2-1)}
             {(\beta_1-\beta_2)^{\alpha_1-\alpha_2-1}\,(\beta_3-\beta_2)^{\alpha_2-\alpha_3-1}\,
             	 \Gamma(2+\alpha_3-\alpha_2)}\,\Phi_1(x)
             \een   
         if $\frac{\pi}{2} < \arg (\beta_3-\beta_1) \leq \frac{3 \pi}{2}$, and  
               \ben
               \Phi^-_{13}(x) &=&
               \Phi^+_{13}(x) \\[0.1ex]
               &+&
               \frac{2 \pi\,i (-1)^{\alpha_1-\alpha_3-3}\,\Gamma(\alpha_1-\alpha_2-1)}
               {(\beta_1-\beta_2)^{\alpha_1-\alpha_2-1}\,(\beta_3-\beta_2)^{\alpha_2-\alpha_3-1}\,
               	\Gamma(2+\alpha_3-\alpha_2)}\,\Phi_1(x)  
               \een   
               if $-\frac{\pi}{2} < \arg (\beta_3-\beta_1) \leq \frac{\pi}{2}$.  
               
           In the same manner one can derive the multiplier $\mu_2$ when
           $\alpha_3-\alpha_1\in\ZZ_{\leq -4}$ but $\alpha_2-\alpha_1\in\ZZ_{\geq -1}$
           or $\alpha_3-\alpha_1\notin\ZZ_{\leq -4}$. Similarly, one can verify that
           $\mu_3$ derived with respect to $\Phi_{13}(x)$ coincides with $\mu_3$
           derived from $\Phi_{23}(x)$.
           
           This completes the proof.    
        \qed

     The Stokes matrices relevant to the actual fundamental matrix constructed by
     \prref{as-3} are the same as these ones given by \thref{Stokes-2} provided that
     $|\beta_3-\beta_1|=|\beta_3-\beta_2|$ and $\R (\alpha_2-\alpha_1) > -1$.
     More precisely,

      \bth{Stokes-3}
  Assume that $\beta_j$'s are distinct such that $|\beta_3-\beta_1| =|\beta_3-\beta_2|$
  but $\beta_3-\beta_1 \neq \pm (\beta_3-\beta_2)$. Assume also that $\R (\alpha_2-\alpha_1) > -1$
  and that  $\arg (\beta_1-\beta_2) \neq \arg (\beta_1-\beta_3) \neq \arg (\beta_2-\beta_3)
  \neq \arg (\beta_1-\beta_2)$.
  Then with
  respect to the actual fundamental matrix $\Phi_{\theta}(x)$ at the origin,
  given by \prref{as-3}, the  equation \eqref{initial} has 
  at most three singular directions $\theta_1=\arg (\beta_1-\beta_2), \theta_2=\arg (\beta_1-\beta_3)$ and
  $\theta_3=\arg (\beta_2-\beta_3)$. The corresponding Stokes matrices are given
  by \thref{Stokes-2} provided that $\Re (\alpha_2-\alpha_1) > -1$
      with the exception of the cases when $\alpha_3-\alpha_1\notin\ZZ_{\leq -4}$
  but $\alpha_3-\alpha_2\in\ZZ_{\geq -1}$ with $\alpha_2-\alpha_1\in\ZZ_{\leq -2}$.
      \ethe


  \section{Iterated integrals as the 1-sum of the product of formal 1-summable series}

    In this section we make a link between the iterated integrals similar to
    $\Phi_{13}(x)$ and summation of the product of two formal 1-summable
     series, that have different singular directions. 

   In \leref{phi} we proved that the formal power series
   $$\,
    \hat{\omega}_1(x)=\sum_{n=0}^{\infty} (-1)^n
    \frac{(2+\alpha_2-\alpha_1)^{(n)}}{(\beta_2-\beta_1)^n}\,x^n\quad
    \textrm{and}\quad
     \hat{\omega}_2(x)=\sum_{n=0}^{\infty} (-1)^n
     \frac{(2+\alpha_3-\alpha_2)^{(n)}}{(\beta_3-\beta_2)^n}\,x^n
   \,$$
  are 1-summable in any directions $\theta_1\neq \arg(\beta_1-\beta_2)$
  and $\theta_2\neq \arg(\beta_2-\beta_3)$, respectively and the functions
   \ben
     \omega_{1, \theta_1}(x)
           &=&
           \int_0^{+\infty e^{i \theta_1}}
     \left(1+\frac{\xi}{\beta_2-\beta_1}\right)^{\alpha_1-\alpha_2-2}
     e^{-\frac{\xi}{x}}\,d \left(\frac{\xi}{x}\right),\\[0,3ex]
        \omega_{2, \theta_2}(x)
        &=&
        \int_0^{+\infty e^{i \theta_2}}
        \left(1+\frac{\xi}{\beta_3-\beta_2}\right)^{\alpha_2-\alpha_3-2}
        e^{-\frac{\xi}{x}}\,d \left(\frac{\xi}{x}\right)
   \een
  define the corresponding 1-sums in such directions.

  Consider the formal series
  \be\label{om}
      \hat{\omega}(x)=\hat{\omega}_1(x)\,\hat{\omega}_2(x).
  \ee
  Note that $\hat{\omega}(x)\in\CC[[x]]_1$ since $\CC[[x]]_1$ is a differential
  sub-algebra of $\CC[[x]]$ stable under product, \cite{L-R}.
  Suppose that $\hat{\omega}(x)$ is  1-summable series. 
  In the differential algebra $\CC\{x\}_1$  the 1-sum of the product 
  $\hat{\omega}_1(x)\,\hat{\omega}_2(x),\,\hat{\omega}_i\in\CC\{x\}_1$ is determined as
  the product of the 1-sums of $\hat{\omega}_1$ and $\hat{\omega}_2$. But what happens when
  $\hat{\omega}_1(x)$ and $\hat{\omega}_2(x)$ have different singular directions? 
  
    Let $\arg(\beta_1-\beta_2) \neq \arg(\beta_2-\beta_3)$. 
   Then following the program  
  of section 4 we can find the 1-sum of $\hat{\omega}(x)$. 
  The formal Borel transform of the series $\hat{\omega}(x)$ is defined as
   the convolution product of the formal Borel transforms 
   $(\hat{\B}_1\,\hat{\omega}_1)(\xi)=\omega_1(\xi)$ and
      $(\hat{\B}_1\,\hat{\omega}_2)(\xi)=\omega_2(\xi)$
      $$\,
        \omega(\xi)=(\hat{\B}_1(\hat{\omega}_1\,\hat{\omega}_2))(\xi)=
        (\omega_1 \star \omega_2)(\xi).
      \,$$
   Moreover, if $R=\min \{|\beta_3-\beta_2|,\,|\beta_2-\beta_1|\}$ then
      \be\label{btr}
      \omega(\xi)=
      (\omega_1 \star \omega_2)(\xi)
          &=&
          \int_0^{\xi}
      \omega_1(z)\,\omega_2(\xi-z)\,d z\nonumber\\[0.15ex]
          &=&
      \int_0^\xi
      \left(1+\frac{z}{\beta_2-\beta_1}\right)^{\alpha_1-\alpha_2-2}
      \left(1+\frac{\xi-z}{\beta_3-\beta_2}\right)^{\alpha_2-\alpha_3-2}\,d z
      \ee
      is a holomorphic in $\{\xi\in\CC\,|\, |\xi| < R\}$ function (see \cite{Sa} for details).
      Once computing $\omega(\xi)$ the corresponding Laplace transform
      yields the 1-sum of the series $\hat{\omega}(x)$ in every non-singular direction
      $\theta$, i.e.
       $$\,
        \omega_{\theta}(x)=\int_0^{+\infty e^{i \theta}}
         \omega(\xi)\,e^{-\frac{\xi}{x}}\,d \left(\frac{\xi}{x}\right).
       \,$$
     This is one of the classical methods of building the 1-sum of the series $\hat{\omega}(x)$.
     But this program is too complicated since the  explicit computation of $\omega(\xi)$ in
     terms of classical functions is out of reach. We propose
     a different approach based on the iterated integrals. Unfortunately
     our approach works only for formal power series $\hat{\omega}_j(x)$ which come from some
     integrals. The asset of this approach, in the case when it works, is that the 1-sum
     obtained by it is convenient  for  measuring  of the Stokes phenomenon. 
     It also detects efficiently the creating of new singular directions in an effect of
     the Borel transform \eqref{btr} (see \cite{Sa-1} for details ).

        \bth{sum}	
        Assume that $\beta_j$'s ate distinct such that $|\beta_3-\beta_1| < |\beta_3-\beta_2|$
        and $|\beta_3-\beta_1| < |\beta_2-\beta_1|$. Assume also that neither
        $\alpha_3-\alpha_1\in\ZZ_{\leq -4}$ nor $\alpha_2-\alpha_1, \alpha_3-\alpha_2\in\ZZ$.
        Then for any directions $\theta_1\neq \arg(\beta_1-\beta_2),\,\theta_2\neq \arg(\beta_1-\beta_3)$
        and $\theta_3 \neq \arg(\beta_2-\beta_3)$ from $0$ to $+\infty \,e^{i \theta_i}$ the function
        \ben
        & &
        \omega_{\theta}(x) =
        - \Gamma(4+\alpha_3-\alpha_1)\,G(\alpha, \beta)
         \int_0^{+\infty e^{i \theta_2}}
        \left(1+\frac{\xi}{\beta_3-\beta_1}\right)^{\alpha_1-\alpha_3-4} e^{-\frac{\xi}{x}}
        d \left(\frac{\xi}{x}\right)\\[0.4ex]
        &+&
        \sum_{s=0}^{\infty}\frac{(2+\alpha_2-\alpha_1)^{(s)}}{(\alpha_2-\alpha_3-1)^{(s)}}
        \left(\frac{\beta_3-\beta_2}{\beta_1-\beta_2}\right)^s
        \int_0^{+\infty e^{i \theta_3}}
        \left(1+\frac{\xi}{\beta_3-\beta_2}\right)^{\alpha_2-\alpha_3-2+s}
        e^{-\frac{\xi}{x}} d \left(\frac{\xi}{x}\right)\\[0.4ex]
        &+&
        \sum_{s=0}^{\infty}\frac{(2+\alpha_3-\alpha_2)^{(s)}}{(\alpha_1-\alpha_2-1)^{(s)}}
        \left(\frac{\beta_2-\beta_1}{\beta_2-\beta_3}\right)^s
        \int_0^{+\infty e^{i \theta_1}}
        \left(1+\frac{\xi}{\beta_2-\beta_1}\right)^{\alpha_1-\alpha_2-2+s}
        e^{-\frac{\xi}{x}} d \left(\frac{\xi}{x}\right),
        \een
        where
        \ben
         G(\alpha, \beta)
            &=& 
             \frac{\Gamma(\alpha_2-\alpha_3-1)}{\Gamma(2+\alpha_2-\alpha_1)}
         \left(\frac{\beta_1-\beta_2}{\beta_3-\beta_2}\right)^{\alpha_1-\alpha_2-2}
         \left(\frac{\beta_1-\beta_3}{\beta_1-\beta_2}\right)^{\alpha_1-\alpha_3-4}
         \\[0.35ex]
           &+&
         \frac{\Gamma(\alpha_1-\alpha_2-1)}{\Gamma(2+\alpha_3-\alpha_2)}
         \left(\frac{\beta_2-\beta_3}{\beta_2-\beta_1}\right)^{\alpha_1-\alpha_2-2}
         \left(\frac{\beta_1-\beta_3}{\beta_2-\beta_3}\right)^{\alpha_1-\alpha_3-4}
        \een
          defines the 1-sum of the series $\hat{\omega}(x)$ in such directions.
        \ethe

        \proof
        
        From \prref{as-2} it follows that when $\alpha_3-\alpha_1\notin\ZZ_{\leq -4}$
        and $\alpha_3-\alpha_2\notin\ZZ$ the actual function $\Phi_{13}(x)$ has the form
        \ben
         & &
        \Phi_{13}(x)\\[0.1ex]
         &=&
        \frac{x^{\alpha_3+4} e^{-\frac{\beta_3}{x}}}{(\beta_3-\beta_2)(\beta_3-\beta_1)}
        \left[F(\alpha; \beta; 0)\,
        \int_0^{+\infty e^{i \theta_2}}
        \left(1+\frac{\nu}{\beta_3-\beta_1}\right)^{\alpha_1-\alpha_3-4}
        e^{-\frac{\nu}{x}} d \left(\frac{\nu}{x}\right)
        - \psi_{\theta}(x)\right],  
        \een
        where $\psi_{\theta}(x)$ is defined by \leref{psi1}. The number series $F(\alpha; \beta; 0)$ 
        is defined by \eqref{F0} for $\sigma=\infty$.
        
        On the other hand applying an integration by parts we find that the formal $\Phi_{13}(x)$
        can be represented as
        \ben
        \Phi_{13}(x)=
        x^{\alpha_1} e^{-\frac{\beta_1}{x}}\left[
        \left(\int_0^x t^{\alpha_2-\alpha_1} e^{-\frac{\beta_2-\beta_1}{t}}\, d t\right)
        \left(\int_0^x t^{\alpha_3-\alpha_2} e^{-\frac{\beta_3-\beta_2}{t}}\, d t\right)
        -\frac{1}{\beta_2-\beta_1}\tilde{I}
        \right],
        \een
        where the first integral is taken in the direction $\arg (\beta_2-\beta_1)$ and the second 
        is taken in the direction $\arg (\beta_3-\beta_2)$. The integral $\tilde{I}$ has the form
        $$\,
        \tilde{I}=\int_0^x t^{\alpha_3-\alpha_1+2}\,e^{-\frac{\beta_3-\beta_1}{t}}
        \left(\sum_{s=0}^{\infty}
        (-1)^s \frac{(2+\alpha_2-\alpha_1)^{(s)}}{(\beta_2-\beta_1)}\,t^s\right)\,d t,
        \,$$
        where the integral is taken in the direction $\arg (\beta_3-\beta_1)$.
        Repeat the methods of \leref{l3} and \leref{psi2} for the intgral $\tilde{I}$
        we obtain a new actual representation of $\Phi_{13}(x)$ provided that
        $|\beta_3-\beta_1| < |\beta_2-\beta_1|$ and $\alpha_2-\alpha_1\notin\ZZ_{\geq -1}$
        \ben
        \Phi_{13}(x) 
        &=&
        \frac{x^{\alpha_3+4} e^{-\frac{\beta_3}{x}}}
        {(\beta_3-\beta_2) (\beta_2-\beta_1)}
        \Big[\omega_{\theta}(x)
        + \frac{\beta_3-\beta_2}{\beta_3-\beta_1} \varphi_{\theta}(x)\\[0.4ex]
        &-&
        \frac{\beta_3-\beta_2}{\beta_3-\beta_1} W(\alpha; \beta; 0)
        \int_0^{+\infty e^{i \theta_2}}
        \left(1+\frac{\nu}{\beta_3-\beta_1}\right)^{\alpha_1-\alpha_3-4}
        e^{-\frac{\nu}{x}} d \left(\frac{\nu}{x}\right)\Big],
        \een
        where $\omega_{\theta}(x)$ is the wanted 1-sum  of the product
        \eqref{om}. Here the infinite number series $W(\alpha; \beta; 0)$ is given by
        $$\,
        W(\alpha; \beta; 0)=\sum_{p=0}^{\infty}
        \frac{(2+\alpha_2-\alpha_1)^{(p)}}{(4+\alpha_3-\alpha_1)^{(p)}}
        \left(\frac{\beta_3-\beta_1}{\beta_2-\beta_1}\right)^p,
        \,$$
        and $\varphi_{\theta}(x)$ is defined by
        \ben
           & &
        \varphi_{\theta}(x) =
        - \frac{\beta_3-\beta_1}{\beta_2-\beta_1}
        \frac{\Gamma(4+\alpha_3-\alpha_1) \Gamma(\alpha_1-\alpha_2-1)}
        {\Gamma(2+\alpha_3-\alpha_2)}
        \left(\frac{\beta_2-\beta_3}{\beta_2-\beta_1}\right)^{\alpha_1-\alpha_2-3}\\[0.4ex]
        &\times&
        \left(\frac{\beta_1-\beta_3}{\beta_2-\beta_3}\right)^{\alpha_1-\alpha_3-4}
        \int_0^{+\infty e^{i \theta_2}}
        \left(1+\frac{\xi}{\beta_3-\beta_1}\right)^{\alpha_1-\alpha_3-4}
        e^{-\frac{\xi}{x}} d\left(\frac{\xi}{x}\right)\\[0.85ex]  
        &+&
        \frac{\beta_3-\beta_1}{\beta_2-\beta_3}
        \sum_{s=0}^{\infty}
        \frac{(2+\alpha_3-\alpha_2)^{(s)}}{(\alpha_1-\alpha_2-1)^{(s)}}
        \left(\frac{\beta_2-\beta_1}{\beta_2-\beta_3}\right)^s
        \int_0^{+\infty\,e^{i \theta_1}}
        \left(1+\frac{\xi}{\beta_2-\beta_1}\right)^{\alpha+s}\,
        e^{-\frac{\xi}{x}}\,d \left(\frac{\xi}{x}\right),
        \een   
        where $\alpha=\alpha_1-\alpha_2-2$.
        Regard the series $F(\alpha; \beta; 0)$
        as the hypergeometric series
        $$\,
        F(a, b; c; z)=\sum_{k=0}^{\infty}
        \frac{(a)^{(k)}\,(b)^{(k)}}{(c)^{(k)}\,k!}\,z^k
        \,$$
        with the parameters $a=1,\,b=2+\alpha_3-\alpha_2,\,c=4+\alpha_3-\alpha_1$
        and $z=\frac{\beta_3-\beta_1}{\beta_3-\beta_2}$. Then from the formula
        $$\,
        F(a, b; c; z) = (1-z)^{-a}\,F\left(a, c-b; c; \frac{z}{z-1}\right)
        \,$$
        it follows that
        $$\,
        F(\alpha; \beta; 0)=\frac{\beta_3-\beta_2}{\beta_1-\beta_2}\,W(\alpha; \beta; 0).
        \,$$     
        Finally, comparing the both actual representation of $\Phi_{13}(x)$
        we arrive at the wanted 1-sum $\omega_{\theta}(x)$ of the formal
        series \eqref{om}.
        
        This completes the proof.
        \qed
        
      In fact \thref{sum} shows that under some assumptions the formal power series
      $\hat{\omega}(x)=\hat{\omega}_1(x)\,\hat{\omega}_2(x)$ is a resurgent formal
      series of Gevrey order 1 \cite{LR-R, Sa-1} in 
      $\CC\backslash\{\beta_1-\beta_2, \beta_2-\beta_3, \beta_1-\beta_3=\beta_1-\beta_2
      +\beta_2-\beta_3\}$ not in   $\CC\backslash\{\beta_1-\beta_2, \beta_2-\beta_3\}$.       
      Then we can say that the direction $\theta=\arg(\beta_1-\beta_3)$ appears as a
      singular direction for the series $\sum_{n=0}^{\infty} a_n x^n$ with $a_n$ from
      \eqref{an} because this series governs a convolution product of two Borel transforms.

      We have a similar result relevant to the case when $\alpha_3-\alpha_1\in\ZZ_{\leq -4}$.
       \bth{sum-1}	
       Assume that  $\beta_j$'s ate distinct such that $|\beta_3-\beta_1| < |\beta_3-\beta_2|$
       and $|\beta_3-\beta_1| < |\beta_2-\beta_1|$. Assume also that 
       $\alpha_3-\alpha_1\in\ZZ_{\leq -4}$ but  $\alpha_2-\alpha_1, \alpha_3-\alpha_2\notin\ZZ_{\leq -2}$.
       Then for any directions $\theta_1\notin \{\arg(\beta_2-\beta_3),\,\arg(\beta_1-\beta_3)\}$ and
       $\theta_2\notin\{ \arg(\beta_1-\beta_2),\,\arg(\beta_1-\beta_3)\}$
        from $0$ to $+\infty \,e^{i \theta_i}$ the function
       \ben
       & &
       (\beta_3-\beta_1)\,\omega_{\theta}(x)\\[0.1ex]
       &=&
        \sum_{l=0}^{\alpha_1-\alpha_3-4}
       [(\beta_2-\beta_1)\,b_l + (\beta_3-\beta_2)\,c_l]\,x^l
       \left(\sum_{s=0}^{\alpha_1-\alpha_3-4-l} (-1)^s
       \frac{(4+l+\alpha_3-\alpha_1)^{(s)}}{(\beta_3-\beta_1)^s}\,x^s\right)
       \\[0.4ex]
       &+&
       x^{\alpha}\Big[
       b_{\alpha} (\beta_2-\beta_1)
       \int_0^{+e^{i \theta_1}}
       \left(1+\frac{\xi}{\beta_3-\beta_2}\right)^{1+\alpha_2-\alpha_1}
       \left(1+\frac{\xi}{\beta_3-\beta_1}\right)^{-1}
        e^{-\frac{\xi}{x}}
       d \left(\frac{\xi}{x}\right)\\[0.4ex]
       &+&
           c_{\alpha} (\beta_3-\beta_2)
           \int_0^{+e^{i \theta_2}}
           \left(1+\frac{\xi}{\beta_2-\beta_1}\right)^{1+\alpha_3-\alpha_2}
           \left(1+\frac{\xi}{\beta_3-\beta_1}\right)^{-1}
           e^{-\frac{\xi}{x}}
           d \left(\frac{\xi}{x}\right)\Big],
           \een
           where $\alpha=\alpha_1-\alpha_3-3$,
       defines the 1-sum of the series $\hat{\omega}(x)$.
       \ethe

      \vspace{1cm}
      {\bf Acknowledgments.}	
      The author was partially supported by Grant  DN 02-5/2016 
      of the Bulgarian Fond ``Scientific Research''.

\begin{small}
    
\end{small}

\end{document}